\documentclass[english,leqno,10pt]{article}
\usepackage{geometry}
\usepackage[T1]{fontenc}
\usepackage[latin9]{inputenc}
\usepackage{babel}
\usepackage{float}
\usepackage{amsmath}
\usepackage{stmaryrd}
\usepackage{amsthm}
\usepackage{amssymb}
\usepackage{cancel}
\usepackage{graphicx,xcolor}
\usepackage{esint}
\usepackage[toc,page]{appendix}
\usepackage{chngcntr}
\usepackage{apptools}
\usepackage{algorithm,algpseudocode}
\newcommand{\R}{{\mathbb R}}
\newcommand{\E}{{\mathbb E}}\makeatletter
\newcommand{\N}{{\mathbb N}}

\floatstyle{ruled}
\newfloat{algorithm}{tbp}{loa}
\providecommand{\algorithmname}{Algorithm}
\floatname{algorithm}{\protect\algorithmname}
\theoremstyle{definition}
\newtheorem*{condition*}{\protect\conditionname}
\theoremstyle{plain}
\newtheorem{thm}{\protect\theoremname}[section]
\theoremstyle{plain}
\newtheorem{prop}[thm]{\protect\propositionname}
\theoremstyle{definition}

\theoremstyle{plain}
\newtheorem{lem}[thm]{\protect\lemmaname}
\theoremstyle{plain}

\theoremstyle{plain}
\newtheorem{remark}[thm]{\protect\remarkname}
\theoremstyle{plain}
\newtheorem{exple}[thm]{\protect\examplename}
\theoremstyle{plain}
\newtheorem*{assumption*}{\protect\assumptionname}
\numberwithin{equation}{section}
\makeatother
\providecommand{\assumptionname}{Assumption}
\providecommand{\conditionname}{Condition}
\providecommand{\corollaryname}{Corollary}
\providecommand{\definitionname}{Definition}
\providecommand{\lemmaname}{Lemma}
\providecommand{\propositionname}{Proposition}
\providecommand{\theoremname}{Theorem}
\providecommand{\remarkname}{Remark}
\providecommand{\examplename}{Example}
\begin{document}

\title{Approximation rate in Wasserstein distance of probability measures on the real line by deterministic empirical measures}
\author{O. Bencheikh and B. Jourdain\thanks{Cermics, \'Ecole des Ponts, INRIA, Marne-la-Vall\'ee, France. E-mails : benjamin.jourdain@enpc.fr, oumaima.bencheikh@enpc.fr. The authors would like to acknowledge financial support from Université Mohammed VI Polytechnique. }}
\maketitle

\begin{abstract}
  We are interested in the approximation in Wasserstein distance with index $\rho\ge 1$ of a probability measure $\mu$ on the real line with finite moment of order $\rho$ by the empirical measure of $N$ deterministic points. The minimal error converges to $0$ as $N\to+\infty$ and we try to characterize the order associated with this convergence. In \cite{xuberger}, Xu and Berger show that, apart when $\mu$ is a Dirac mass and the error vanishes, the order is not larger than $1$ and give a sufficient condition for the order to be equal to this threshold $1$ in terms of the density of the absolutely continuous with respect to the Lebesgue measure part of $\mu$. They also prove that the order is not smaller than $1/\rho$ when the support of $\mu$ is bounded and not larger when the support is not an interval. We complement these results by checking that for the order to lie in the interval $\left(1/\rho,1\right)$, the support has to be bounded and by stating a necessary and sufficient condition in terms of the tails of $\mu$ for the order to be equal to some given value in the interval $\left(0,1/\rho\right)$, thus precising the sufficient condition in terms of moments given in \cite{xuberger}. In view of practical application, we emphasize that in the proof of each result about the order of convergence of the minimal error, we exhibit a choice of points explicit in terms of the quantile function of $\mu$ which exhibits the same order of convergence. 

\bigskip
  \noindent{\bf Keywords:}
  deterministic empirical measures, Wasserstein distance, rate of convergence.

\bigskip
\noindent {{\bf AMS Subject Classification (2010):} \it 49Q22, 60-08}
\end{abstract}


\section*{Introduction}

Let $\rho\ge 1$ and $\mu$ be a probability measure on the real line. We are interested in the rate of convergence in terms of $N\in\N^*$ of 
\begin{equation}
  e_N(\mu,\rho):=\inf\left\{\mathcal{W}_\rho\left(\frac{1}{N}\sum_{i=1}^N\delta_{x_i},\mu\right):-\infty<x_1\le x_2\le \cdots\le x_N<+\infty\right\},\label{pbor}
\end{equation}
where ${\cal W}_\rho$ denotes the Wasserstein distance with index $\rho$. The motivation is the approximation of the probability measure $\mu$ by finitely supported probability measures. An example of application is provided by the optimal initialization of systems of particles with mean-field interaction \cite{jrdcds,benchjour}, where, to preserve the mean-field feature, it is important to get $N$ points with equal weight $\frac{1}{N}$ (of course, nothing prevents several of these points to be equal).
The Hoeffding-Fr\'echet or comonotone coupling between two probability measures $\nu$ and $\eta$ on the real line is optimal for $\mathcal{W}_\rho$ so that:
\begin{align}\label{Wasserstein}
  \displaystyle \mathcal{W}_\rho^\rho\left(\nu,\eta\right) = \int_0^1 \left|F^{-1}_{\nu}(u) - F^{-1}_{\eta}(u) \right|^\rho\,du,
\end{align}
where for $u\in(0,1)$, $F^{-1}_{\nu}(u)= \inf\left\{ x \in \R: \nu\left(\left(-\infty,x\right]\right)\geq u\right\}$ and $F^{-1}_{\eta}(u)= \inf\left\{ x \in \R: \eta\left(\left(-\infty,x\right]\right)\geq u\right\}$ are the respective quantile functions of $\nu$ and $\eta$. We set $F(x)=\mu\left(\left(-\infty,x\right]\right)$ for $x\in\R$ and denote $F^{-1}(u)=\inf\left\{ x \in \R: F(x)\geq u\right\}$ for $u\in(0,1)$. We have $u\le F(x)\Leftrightarrow F^{-1}(u)\le x$. The quantile function $F^{-1}$ is left-continuous and non-decreasing and we denote by $F^{-1}(u+)$ its right-hand limit at $u\in [0,1)$ (in particular $F^{-1}(0+)=\lim\limits_{u\to 0+}F^{-1}(u)\in [-\infty,+\infty)$) and set $F^{-1}(1)=\lim\limits_{u\to 1-}F^{-1}(u)\in(-\infty,+\infty]$.

By \eqref{Wasserstein}, when $-\infty<x_1\le x_2\le \cdots\le x_N<+\infty$,   
\begin{equation}
	\mathcal{W}_\rho^\rho\left(\frac{1}{N}\sum_{i=1}^N\delta_{x_i},\mu\right) =\sum \limits_{i=1}^N \int_{\frac{i-1}{N}}^{\frac{i}{N}} \left|x_{i}-F^{-1}(u)\right|^\rho\,du,\label{wrhomunmu}
\end{equation}
where, by the inverse transform sampling, the right-hand side is finite if and only if $\int_\R|x|^\rho\mu(dx)<+\infty$. So, when considering $e_N(\mu,\rho)$, we will suppose that $\mu$ has a finite moment of order $\rho$.\\

In the first section of the paper, we recall that, under this moment condition, the infimum in \eqref{pbor} is attained: 
\begin{equation*}
   e_N\left(\mu,\rho\right) = {\cal W}_\rho\left(\frac{1}{N}\sum_{i=1}^N\delta_{x_i^N},\mu\right)=\sum \limits_{i=1}^N \int_{\frac{i-1}{N}}^{\frac{i}{N}} \left|x^N_{i}-F^{-1}(u)\right|^\rho\,du
\end{equation*} 
for some points $x_i^N\in[F^{-1}(\frac{i-1}{N}+),F^{-1}(\frac{i}{N})]\cap\R$ which are unique as soon as $\rho>1$ and explicit in the quadratic case $\rho=2$. Of course the points $\left(x_i^N\right)_{1\le i\le N}$ depend on $\rho$ but we do not explicit this dependence to keep notations simple. For $\rho=1$, because of the lack of strict convexity of $\R\ni x\mapsto|x|$, there may be several optimal choices among which $x_i^N=F^{-1}\left(\frac{2i-1}{2N}\right)$ for $i\in\{1,\hdots,N\}$.  Note that when $\rho\ge\tilde \rho\ge 1 $ and $\displaystyle \int_\R|x|^\rho\mu(dx)<+\infty$, with $(x_i^N)_{1\le i\le N}$ denoting the optimal points for $\rho \ge 1$,
\begin{equation}
   e_N(\mu,\rho)= {\cal W}_\rho\left(\frac{1}{N}\sum_{i=1}^N\delta_{x^N_i},\mu\right)\ge {\cal W}_{\tilde \rho}\left(\frac{1}{N}\sum_{i=1}^N\delta_{x^N_i},\mu\right)\ge e_N(\mu,\tilde\rho)\label{minoerhoun}.
 \end{equation}Hence $\rho\mapsto e_N(\mu,\rho)$ is non-decreasing.

We give an alternative expression of $e_N(\mu,\rho)$ in terms of the cumulative distribution function rather than the quantile function and recover that $e_N(\mu,\rho)$ tends to $0$ as $N\to+\infty$ when $\int_\R|x|^\rho\mu(dx)<+\infty$. The main purpose of the paper is to study the rate at which this convergence occurs. In particular, we would like to give sufficient conditions on $\mu$, which, when possible, are also necessary, to ensure convergence at a rate $N^{-\alpha}$ with $\alpha >0$ called the order of convergence.


This question has already been studied by Xu and Berger in \cite{xuberger}. According to Theorem 5.20 \cite{xuberger},  $\limsup_{N\to\infty}Ne_N(\mu,\rho)\ge\frac 12\int_0^{\frac 12}|F^{-1}(u+\frac 12)-F^{-1}(u)|du$ with the right-hand side positive apart when $\mu$ is a Dirac mass and $e_N(\mu,\rho)$ vanishes for all $N,\rho\ge 1$. This result may be complemented by the non-asymptotic bound obtained in Lemma 1.4 of the first version \cite{BJ} of the present paper : $Ne_N(\mu,\rho)+(N+1)e_{N+1}(\mu,\rho)\ge \frac 12\int_{\R}F(x)\wedge (1-F(x))\,dx$ where $\int_{\R}F(x)\wedge (1-F(x))\,dx=\int_0^{\frac 12}|F^{-1}(u+\frac 12)-F^{-1}(u)|du$ as easily seen since these integrals correspond to the area of the points at the right to $(F(x))_{-\infty\le x\le F^{-1}(\frac{1}{2})}$, at the left to $(1-F(x))_{F^{-1}(\frac{1}{2})<x<+\infty}$, above $0$ and below $\frac{1}{2}$ respectively computed by integration with respect to the abscissa and to the ordinate.
This previous version also contained a section devoted to the case when the support of $\mu$ is bounded i.e. $F^{-1}(1)-F^{-1}(0+)<+\infty$. The results in this section were mainly obtained before by Xu and Berger in \cite{xuberger}. In particular, according to Theorem 5.21 (ii) \cite{xuberger}, when the support of $\mu$ is bounded, then $\sup_{N\ge 1}N^{\frac{1}{\rho}}e_N(\mu,\rho)<+\infty$, while when $F^{-1}$ is discontinuous $\limsup_{N\to+\infty}N^{\frac{1}{\rho}}e_N(\mu,\rho)>0$, according to Remark 5.22 (ii) \cite{xuberger}. We recovered these results and went slightly further by noticing that when the support of $\mu$ is bounded and $F^{-1}$ is continuous, then $\lim_{N\to\infty}N^{\frac{1}{\rho}}e_N(\mu,\rho)=0$ (Proposition 2.1 \cite{BJ}), with an order of convergence arbitrarily slow as examplified by the the beta distribution with parameter $(\beta,1)$ with $\beta>0$ : $\mu_\beta(dx)=\beta {\mathbf 1}_{[0,1]}(x)x^{\beta-1}\,dx$. Indeed, with the notation $\asymp$ defined at the end of the introduction, according to Example 2.3 \cite{BJ}, when $\rho>1$, $e_N(\mu_\beta,\rho)\asymp N^{-\frac{1}{\rho}-\frac{1}{\beta}}\asymp{\cal W}_\rho\left(\frac{1}{N}\sum_{i=1}^N\delta_{F^{-1}\left(\frac{i-1}{N}\right)},\mu_\beta\right)$ for $\beta>\frac{\rho}{\rho-1}$, $e_N(\mu_\beta,\rho)\asymp N^{-1}(\ln N)^{\frac{1}{\rho}}\asymp{\cal W}_\rho\left(\frac{1}{N}\sum_{i=1}^N\delta_{F^{-1}\left(\frac{i-1}{N}\right)},\mu_\beta\right)$ for $\beta=\frac{\rho}{\rho-1}$ and $\lim_{N\to+\infty}Ne_N(\mu,\rho)=\frac{1}{2\beta(\rho+1)^{1/\rho}}\left(\int_0^1u^{\frac{\rho}{\beta}-\rho}du\right)^{1/\rho}$ when $\beta\in (0,\frac{\rho}{\rho-1})$. The latter limiting behaviour which remains valid for $\rho=1$ whatever $\beta>0$ is a consequence of one of the main results by Xu and Berger \cite{xuberger} namely Theorem 5.15 : when the density $f$ of the absolutely continuous with respect to the Lebesgue measure part of $\mu$ is $dx$ a.e. positive on $\left\{x\in\R:0<F(x)<1\right\}$ (or equivalently $F^{-1}$ is absolutely continuous), then
$$\lim_{N\to+\infty} Ne_N(\mu,\rho)=\lim_{N\to+\infty} N{\cal W}_\rho\left(\frac{1}{N}\sum_{i=1}^N\delta_{F^{-1}\left(\frac{2i-1}{2N}\right)},\mu\right) = \frac{1}{2(\rho+1)^{1/\rho}}\left(\int_\R\frac{{\mathbf 1}_{\{0<F(x)<1\}}}{f^{\rho-1}(x)}\,dx\right)^{1/\rho}.$$
In Theorem 2.4 \cite{BJ}, we recovered this result and also stated that, without the positivity assumption on the density, 
$$\liminf_{N\to+\infty} Ne_N(\mu,\rho) \ge \frac{1}{2(\rho+1)^{1/\rho}} \left(\int_\R\frac{{\mathbf 1}_{\{f(x)>0\}}}{f^{\rho-1}(x)}\,dx\right)^{1/\rho}.$$
In particular for $(Ne_N(\mu,1))_{N\ge 1}$ to be bounded the Lebesgue measure of $\{x\in\R:f(x)>0\}$ must be finite. Weakening this necessary condition is discussed in Section 4 of the first version \cite{BJ} of this paper. In \cite{chevallier}, Chevallier addresses the multidimensional setting and proves in Theorem III.3 that 
for a probability measure $\mu$ on $\R^d$ with support bounded by $r$, there exist points $x_1,\hdots,x_N\in\R^d$ such that $\frac{1}{4r}\mathcal{W}_\rho\left(\frac{1}{N}\sum_{i=1}^N\delta_{x_i},\mu\right)\le f_{\rho,d}(N)$ where $f_{\rho,d}(N)$ is respectively equal to $\left(\frac d{d-\rho}\right)^{\frac 1\rho}N^{-\frac 1 d}$, $\left(\frac{1+\ln N}{N}\right)^{\frac 1 d}$, and $\zeta(p/d)N^{-\frac 1 \rho}$ with $\zeta$ denoting the zeta Riemann function when $\rho<d$, $\rho=d$ and $\rho>d$.

The case when the support of $\mu$ is not bounded is also considered by Xu and Berger \cite{xuberger} in the one-dimensional setting of the present paper and by  \cite{chevallier} in the multidimensional setting. In Corollary III.5 \cite{chevallier}, Chevallier proves that $\lim_{N\to\infty}(f_{\rho,d}(N))^{-\alpha\rho}\mathcal{W}_\rho\left(\frac{1}{N}\sum_{i=1}^N\delta_{x_i},\mu\right)=0$ when $\int_{\R^d}|x|^{\frac{\rho}{1-\alpha\rho}}\mu(dx)<+\infty$ for some $\alpha\in (0,\frac 1 \rho)$. This generalizes the one-dimensional statement of Theorem 5.21 (i) \cite{xuberger} : under the same moment condition, $\lim_{N\to\infty}N^\alpha e_N(\mu,\rho)=0$.   In Theorem \ref{alphaRater}, using our alternative formula for $e_N(\mu,\rho)$ in terms of the cumulative distribution function of $\mu$, we refine this result by stating the following necessary and sufficient condition
$$\forall \alpha\in \left(0,\frac{1}{\rho}\right),\;\lim_{x\to +\infty}x^{\frac{\rho}{1-\alpha \rho}}\Big(F(-x)+1-F(x)\Big)=0\Leftrightarrow \lim_{N\to+\infty}N^{\alpha}e_N(\mu,\rho)=0.$$
We also check that
$$\forall \alpha\in \left(0,\frac{1}{\rho}\right),\quad\displaystyle{\sup_{N \ge 1}} N^{\alpha} \, e_N(\mu,\rho)<+\infty \Leftrightarrow \sup_{x\ge 0}\;x^{\frac{\rho}{1-\alpha \rho}}\Big(F(-x)+1-F(x)\Big)<+\infty,$$
a condition under which, the order of convergence $\alpha$ of the minimal error $e_N(\mu,\rho)$ is preserved by choosing $x_1=F^{-1}\left(\frac{1}{N}\right)\wedge (-N^{\frac{1}{\rho}-\alpha})$, $x_N=F^{-1}\left(\frac{N-1}{N}\right)\vee N^{\frac{1}{\rho}-\alpha}$ and any $x_i\in[F^{-1}(\frac{i-1}{N}+),F^{-1}(\frac{i}{N})]$ for $2\le i\le N-1$.
 We also exhibit non-compactly supported probability measures $\mu$ such that, for $\rho>1$, $\lim_{N\to+\infty}N^{\frac{1}{\rho}} \, e_N(\mu,\rho)=0$. Nevertheless, we show that for $\left(N^\alpha e_N(\mu,\rho)\right)_{N\ge 1}$ to be bounded for $\alpha>\frac{1}{\rho}$, the support of $\mu$ has to be bounded. We last give a necessary condition for $\left(N^{\frac{1}{\rho}} e_N(\mu,\rho)\right)_{N\ge 1}$ to be bounded, which unfortunately is not sufficient but ensures the boundedness of  $\left(\frac{N^{1/\rho}}{1+\ln N}e_N(\mu,\rho)\right)_{N\ge 1}$.

We summarize our results together with the ones obtained by Xu and Berger \cite{xuberger} in Table \ref{res}.

\begin{center}
\begin{table}[!ht]
\begin{tabular}{|c||c|c|}\hline
  $\alpha$ & Necessary condition & Sufficient condition\\\hline\hline
  $\alpha=1$ & $\displaystyle \int_\R\frac{{\mathbf 1}_{\left\{f(x)>0 \right\}}}{f^{\rho-1}(x)}\,dx<+\infty$ & $f(x)>0$ $dx$ a.e. on $\left\{x\in\R: 0<F(x)<1 \right\}$ \\
   & (Thm. 2.4 \cite{BJ})  &  and $\displaystyle \int_\R\frac{{\mathbf 1}_{\{f(x)>0\}}}{f^{\rho-1}(x)}\,dx<+\infty$ (Thm. 5.15 \cite{xuberger}) \\\hline
  $\alpha\in \left(\frac 1\rho,1\right)$& $F^{-1}$ continuous (Remark 5.22 (ii) \cite{xuberger})  & related to the modulus of continuity of $F^{-1}$\\
          when $\rho>1$     & and $\mu$ with bounded support (Prop. \ref{propals1rcomp}) & \\\hline
  $\alpha=\frac 1\rho$ & $\exists \lambda>0$, $\forall x\ge 0$, $F(-x)+1-F(x)\le\frac{e^{-\lambda x}}{\lambda}$ & $\mu$  with bounded support (Thm. 5.21 (ii) \cite{xuberger}) \\
             & (Prop. \ref{propal1rho}) & For $\rho>1$,  Exple \ref{exempleexp} with unbounded supp.\\\hline
  $\alpha\in\left(0,\frac 1\rho\right)$ &
${\sup\limits_{x\ge 0}}\;x^{\frac{\rho}{1-\alpha \rho}}\Big(F(-x)+1-F(x)\Big)<+\infty$& ${\sup\limits_{x\ge 0}}\;x^{\frac{\rho}{1-\alpha \rho}}\Big(F(-x)+1-F(x)\Big)<+\infty$\\ &(Thm. \ref{alphaRater}) &(Thm. \ref{alphaRater})\\\hline
\end{tabular}\caption{Conditions for the convergence of $e_N(\mu,\rho)$ with order $\alpha$ : $\sup\limits_{N\ge 1}N^\alpha e_N(\mu,\rho)<+\infty$.}\label{res}
\end{table}   
\end{center}
{\bf Notation :} \begin{itemize}
\item We denote by $\lfloor x\rfloor$ (resp. $\lceil x\rceil$) the integer $j$ such that $j\le x<j+1$ (resp. $j-1<x\le j$) and by $\{x\}=x-\lfloor x\rfloor$ the integer part of $x\in\R$.
  \item For two sequences $(a_N)_{N\ge 1}$ and $(b_N)_{N\ge 1}$ of real numbers with $b_N>0$ for $N\ge 2$ we denote $a_N\asymp b_N$ when $\displaystyle 0<\inf_{N\ge 2}\left(\frac{a_N}{b_N}\right)$ and $\displaystyle \sup_{N\ge 2}\left(\frac{a_N}{b_N}\right)<+\infty$.
\end{itemize}


\section{Preliminary results}

When $\rho=1$ (resp. $\rho=2$), $\displaystyle \R\ni y\mapsto N\int_{\frac{i-1}{N}}^{\frac{i}{N}} \left|y-F^{-1}(u)\right|^\rho\,du$ is minimal for $y$ belonging to the set $\left[F^{-1}\left(\frac{2i-1}{2N}\right),F^{-1}\left(\frac{2i-1}{2N}+\right)\right]$ of medians (resp. equal to the mean $\displaystyle N\int_{\frac{i-1}{N}}^{\frac{i}{N}}F^{-1}(u)\,du$) of the image of the uniform law on $\left[\frac{i-1}{N},\frac{i}{N}\right]$ by $F^{-1}$. For general $\rho>1$, the function $\displaystyle \R\ni y\mapsto \int_{\frac{i-1}{N}}^{\frac{i}{N}}\left|y-F^{-1}(u)\right|^\rho\,du$ is strictly convex and continuously differentiable with derivative 
\begin{equation}
   \rho\int_{\frac{i-1}{N}}^{\frac{i}{N}}\left({\mathbf 1}_{\left\{y\ge F^{-1}(u)\right\}} \left(y-F^{-1}(u)\right)^{\rho-1}-{\mathbf 1}_{\left\{y<F^{-1}(u)\right\}}\left(F^{-1}(u)-y\right)^{\rho-1}\right)\,du\label{dery}
\end{equation} 
non-positive for $y=F^{-1}\left(\frac{i-1}{N}+\right)$ when either $i=1$ and $F^{-1}(0+)>-\infty$ or $i\ge 2$ and non-negative for $y=F^{-1}\left(\frac{i}{N}\right)$ when either $i\le N-1$ or $i=N$ and $F^{-1}(1)<+\infty$. Since the derivative has a positive limit as $y\to+\infty$ and a negative limit as $y\to-\infty$, we deduce that $\displaystyle \R\ni y\mapsto \int_{\frac{i-1}{N}}^{\frac{i}{N}}\left|y-F^{-1}(u)\right|^\rho\,du$ admits a unique minimizer $x_i^N\in \left[F^{-1}\left(\frac{i-1}{N}+\right),F^{-1}\left(\frac{i}{N}\right)\right]\cap\R$ (to keep notations simple, we do not explicit the dependence of $x_i^N$ on $\rho$). Therefore 
\begin{equation}
  e^\rho_N(\mu,\rho)=\sum_{i=1}^N \int_{\frac{i-1}{N}}^{\frac{i}{N}} \left|x_{i}^N-F^{-1}(u)\right|^\rho\,du \mbox{ with } \left[F^{-1}\left(\frac{i-1}{N}+\right),F^{-1}\left(\frac{i}{N}\right)\right]\ni x_i^N=\begin{cases}
    \displaystyle F^{-1}\left(\frac{2i-1}{2N}\right)\mbox{ if }\rho=1,\\
    \displaystyle N\int_{\frac{i-1}{N}}^{\frac{i}{N}}F^{-1}(u)\,du\mbox{ if }\rho=2,\\
    \mbox{not explicit otherwise.}
  \end{cases}\label{enrho}
\end{equation}
When needing to bound $e_N(\mu,\rho)$ from above, we may replace the optimal point $x_i^N$ by $F^{-1}\left(\frac{2i-1}{2N}\right)$:
$$\forall i\in\{1,\hdots,N\},\quad \int_{\frac{i-1}{N}}^{\frac{i}{N}}\left|F^{-1}(u)-x_i^N \right|^\rho\,du\le \int_{\frac{i-1}{N}}^{\frac{i}{N}}\left|F^{-1}(u)-F^{-1}\left(\frac{2i-1}{2N}\right)\right|^\rho\,du,$$
a simple choice particularly appropriate when linearization is possible since $\displaystyle \left[\frac{i-1}{N},\frac{i}{N}\right]\ni v\mapsto\int_{\frac{i-1}{N}}^{\frac{i}{N}}\left|u-v\right|^\rho\,du$ is minimal for $v=\frac{2i-1}{2N}$.To bound $e_N(\mu,\rho)$ from below, we can use that, by Jensen's inequality and the minimality of $F^{-1}\left(\frac{2i-1}{2N}\right)$ for $\rho=1$,
\begin{align}
 \displaystyle \int_{\frac{i-1}{N}}^{\frac i N} \left|F^{-1}(u)-x^N_{i}\right|^{\rho}\,du &\ge N^{\rho -1} \left(\int_{\frac{i-1}{N}}^{\frac i N} \left|F^{-1}(u)-x^N_{i}\right|\,du\right)^{\rho} \ge N^{\rho -1} \left(\int_{\frac{i-1}{N}}^{\frac i N} \left|F^{-1}(u)-F^{-1}\left(\frac{2i-1}{2N}\right)\right|\,du\right)^{\rho}  \notag\\
 &\ge N^{\rho -1} \left(\frac{1}{4N}\left(F^{-1}\left(\frac{2i-1}{2N}\right)-F^{-1}\left(\frac{4i-3}{4N}\right)+F^{-1}\left(\frac{4i-1}{4N}\right)-F^{-1}\left(\frac{2i-1}{2N}\right)\right)\right)^{\rho} \notag\\
 &\ge \frac{1}{4^\rho N}\left(F^{-1}\left(\frac{4i-1}{4N}\right)-F^{-1}\left(\frac{4i-3}{4N}\right)\right)^{\rho}.\label{minotermbordlimB}
\end{align}

We also have an alternative formulation of $e_N(\mu,\rho)$ in terms of the cumulative distribution function $F$ in place of the quantile function $F^{-1}$:
\begin{prop}\label{propenf}
  \begin{equation}
    e_N^\rho(\mu,\rho)=\rho\sum_{i=1}^N\left(\int_{F^{-1}\left(\frac{i-1}{N}+\right)}^{x_i^N}\left(x_i^N-y\right)^{\rho-1}\left(F(y)-\frac{i-1}{N}\right)\,dy+\int^{F^{-1}\left(\frac{i}{N}\right)}_{x_i^N}\left(y-x_i^N\right)^{\rho-1}\left(\frac{i}{N}-F(y)\right)\,dy\right).\label{enf}
  \end{equation}
\end{prop} 
Under the convention $F^{-1}(0)=-\infty$, when, for some $i\in\{1,\hdots,N\}$,  $F^{-1}\left(\frac{i-1}{N}+\right)>F^{-1}\left(\frac{i-1}{N}\right)$, then $F(y)=\frac{i-1}{N}$ for $y\in\left[F^{-1}\left(\frac{i-1}{N}\right),F^{-1}\left(\frac{i-1}{N}+\right)\right)$ and $\displaystyle \int_{F^{-1}\left(\frac{i-1}{N}\right)}^{F^{-1}\left(\frac{i-1}{N}+\right)}(x_i^N-y)^{\rho-1}\left(F(y)-\frac{i-1}{N}\right)\,dy=0$ so that the lower integration limit in the first integral in the right-hand side of \eqref{enf} may be replaced by $F^{-1}\left(\frac{i-1}{N}\right)$. In a similar way, the upper integration limit in the second integral may be replaced by $F^{-1}\left(\frac{i}{N}+\right)$ under the convention $F^{-1}(1+)=+\infty$.

When $\rho=1$, the equality \eqref{enf} follows from the interpretation of ${\cal W}_1(\nu,\eta)$ as the integral of the absolute difference between the cumulative distribution functions of $\nu$ and $\eta$ (equal, as seen with a rotation with angle $\frac{\pi}{2}$, to the integral of the absolute difference between their quantile functions) and the integral simplifies into:
\begin{equation}\label{w1altern2b}
  e_N(\mu,1)=\sum_{i=1}^N\left(\int_{F^{-1}\left(\frac{i-1}{N}+\right)}^{F^{-1}\left(\frac{2i-1}{N}\right)}\left(F(y)-\frac{i-1}{N}\right)\,dy+\int^{F^{-1}\left(\frac{i}{N}\right)}_{F^{-1}\left(\frac{2i-1}{N}\right)}\left(\frac{i}{N}-F(y)\right)\,dy\right)=\frac{1}{N}\int_\R\min_{j\in\N}\left|NF(y)-j\right|\,dy.
\end{equation}
For $\rho>1$, it can be deduced from the general formula for ${\cal W}_\rho^\rho(\nu,\eta)$ in terms of the cumulative distribution functions of $\mu$ and $\eta$ (see for instance Lemma B.3 \cite{jourey2}). It is also a consequence of the following equality for each term of the decomposition over $i\in\{1,\hdots,N\}$, that we will need next.
\begin{lem}\label{lemenf}
  Assume that $\displaystyle \int_\R|x|^\rho\mu(dx)<+\infty$ with $\rho\ge 1$. For $i\in\{1,\hdots,N\}$ and $x\in \left[F^{-1}\left(\frac{i-1}{N}\right),F^{-1}\left(\frac{i}{N}\right)\right]\cap\R$ (with convention $F^{-1}(0)=F^{-1}(0+)$), we have:
  $$\int_{\frac{i-1}{N}}^{\frac{i}{N}}\left|x-F^{-1}(u)\right|^\rho\,du=\rho\int_{F^{-1}\left(\frac{i-1}{N}\right)}^{x}(x-y)^{\rho-1}\left(F(y)-\frac{i-1}{N}\right)\,dy+\rho\int^{F^{-1}\left(\frac{i}{N}\right)}_{x}(y-x)^{\rho-1}\left(\frac{i}{N}-F(y)\right)\,dy,$$
  and the right-hand side is minimal for $x=x_i^N$.
\end{lem}
\begin{proof} 
Let $i\in\{1,\hdots,N\}$ and $x\in\left[F^{-1}\left(\frac{i-1}{N}\right),F^{-1}\left(\frac{i}{N}\right)\right]\cap\R$. We have $\frac{i-1}{N}\le F(x)$ and $F(x-)\le\frac{i}{N}$. Since $F^{-1}(u)\le x\Leftrightarrow u\le F(x)$ and $F^{-1}(u)=x$ for $u\in \left(F(x-),F(x)\right]$, we have:
$$\int_{\frac{i-1}{N}}^{\frac{i}{N}} \left|x-F^{-1}(u)\right|^\rho\,du=\int_{\frac{i-1}{N}}^{F(x)} \left(x-F^{-1}(u)\right)^\rho\,du+\int^{\frac{i}{N}}_{F(x)} \left(F^{-1}(u)-x\right)^\rho\,du.$$
Using the well-known fact that the image of ${\mathbf 1}_{[0,1]}(v)\,dv\mu(dz)$ by $(v,z)\mapsto F(z-)+v\mu(\{z\})$ is the Lebesgue measure on $[0,1]$ and that ${\mathbf 1}_{[0,1]}(v)\,dv\mu(dz)$ a.e., $F^{-1}\left(F(z-)+v\mu(\{z\})\right)=z$ , we obtain that:
\begin{align}
  \int_{\frac{i-1}{N}}^{F(x)} \left(x-F^{-1}(u)\right)^\rho\,du &= \int_{v=0}^1\int_{z\in \R}{\mathbf 1}_{\left\{\frac{i-1}{N}\le F(z-)+v\mu(\{z\})\le F(x)\right\}}(x-z)^\rho\mu(dz)\,dv\notag\\
  &=\int_{v=0}^1\int_{z\in \R}{\mathbf 1}_{\left\{\frac{i-1}{N}\le F(z-)+v\mu(\{z\})\le F(x)\right\}}\int \rho(x-y)^{\rho-1}{\mathbf 1}_{\{z\le y\le x\}}\,dy\mu(dz)\,dv\notag\\
  &=\rho\int_{y=-\infty}^{x}(x-y)^{\rho-1}\int_{v=0}^1\int_{z\in \R}{\mathbf 1}_{\left\{\frac{i-1}{N}\le F(z-)+v\mu(\{z\})\right\}}{\mathbf 1}_{\left\{z\le y \right\}}\mu(dz)\,dv\,dy.\label{triplint}
\end{align}
For $v>0$, $\{z\in\R:F(z-)+v\mu(\{z\})\le F(y)\}=(-\infty,y]\cup\{z\in\R:z>y\mbox{ and }F(z)=F(y)\}$ with $\mu\left(\{z\in\R:z>y\mbox{ and }F(z)=F(y)\}\right)=0$ and therefore
$$\int_{z\in \R}{\mathbf 1}_{\left\{\frac{i-1}{N}\le F(z-)+v\mu(\{z\})\right\}}{\mathbf 1}_{\{z\le y\}}\mu(dz)=\int_{z\in \R}{\mathbf 1}_{\left\{\frac{i-1}{N}\le F(z-)+v\mu(\{z\})\le F(y)\right\}}\mu(dz).$$ 
Plugging this equality in \eqref{triplint}, using again the image of ${\mathbf 1}_{[0,1]}(v)\,dv\mu(dz)$ by $(v,z)\mapsto F(z-)+v\mu(\{z\})$ and the equivalence $\frac{i-1}{N}\le F(y)\Leftrightarrow F^{-1}\left(\frac{i-1}{N}\right)\le y$ 
, we deduce that:
\begin{align*}
  \int_{\frac{i-1}{N}}^{F(x)}\left(x-F^{-1}(u)\right)^\rho\,du&=\rho\int_{y=-\infty}^{x}(x-y)^{\rho-1}\int_{u=0}^1{\mathbf 1}_{\left\{\frac{i-1}{N}\le u\le F(y)\right\}}\,du\,dy=\rho\int_{F^{-1}\left(\frac{i-1}{N}\right)}^{x}(x-y)^{\rho-1}\left(F(y)-\frac{i-1}{N}\right)\,dy.
\end{align*}
In a similar way, we check that: 
$$\int^{\frac{i}{N}}_{F(x)}\left(F^{-1}(u)-x\right)^\rho\,du=\rho\int^{F^{-1}\left(\frac{i}{N}\right)}_{x}(y-x)^{\rho-1}\left(\frac{i}{N}-F(y)\right)\,dy,$$
which concludes the proof.\end{proof}
 
\begin{prop}
  For each $\rho\ge 1$, we have $\displaystyle \int_\R |x|^\rho\mu(dx)<+\infty\Leftrightarrow \lim_{N\to+\infty} e_N(\mu,\rho)=0$. 
\end{prop}

The direct implication can be deduced from the inequality $e_N(\mu,\rho)\le{\cal W}_\rho\left(\frac{1}{N}\sum\limits_{i=1}^N\delta_{X_i},\mu\right)$ and the almost sure convergence to $0$ of ${\cal W}_\rho\left(\frac{1}{N}\sum\limits_{i=1}^N\delta_{X_i},\mu\right)$ for $(X_i)_{i\ge 1}$ i.i.d. according to $\mu$ deduced from the strong law of large numbers and stated for instance in Theorem 2.13 \cite{bobkovledoux}. We give an alternative simple argument based on \eqref{enf}.
\begin{proof}
According to the introduction, the finiteness of $e_N(\mu,\rho)$ for some $N\ge 1$ implies that $\displaystyle \int_\R|x|^\rho\mu(dx)<+\infty$. So it is enough to check the zero limit property under the finite moment condition. 

When respectively $F^{-1}\left(\frac{i}{N}\right)\le 0$, $F^{-1}\left(\frac{i-1}{N}+\right)<0<F^{-1}\left(\frac{i}{N}\right)$ or $F^{-1}\left(\frac{i-1}{N}+\right)\ge 0$ , then, by Lemma \ref{lemenf}, the term with index $i$ in \eqref{enf} is respectively bounded from above by 
$$\int_{F^{-1}\left(\frac{i-1}{N}+\right)}^{F^{-1}\left(\frac{i}{N}\right)}\left(F^{-1}\left(\frac{i}{N}\right)-y\right)^{\rho-1}\left(F(y)-\frac{i-1}{N}\right)\,dy \le \int_{F^{-1}\left(\frac{i-1}{N}+\right)}^{F^{-1}\left(\frac{i}{N}\right)}(-y)^{\rho-1}\left(\frac{1}{N}\wedge F(y)\right)\,dy,$$
$$\int_{F^{-1}\left(\frac{i-1}{N}+\right)}^{0}(-y)^{\rho-1}\left(\frac{1}{N}\wedge F(y)\right)\,dy+\int^{F^{-1}\left(\frac{i}{N}\right)}_{0}y^{\rho-1}\left(\frac{1}{N}\wedge(1-F(y))\right)\,dy,$$
$$\int_{F^{-1}\left(\frac{i-1}{N}+\right)}^{F^{-1}\left(\frac{i}{N}\right)}\left(y-F^{-1}\left(\frac{i-1}{N}+\right)\right)^{\rho-1}\left(\frac{i}{N}-F(y)\right)\,dy \le \int_{F^{-1}\left(\frac{i-1}{N}+\right)}^{F^{-1}\left(\frac{i}{N}\right)}y^{\rho-1}\left(\frac{1}{N}\wedge(1-F(y))\right)\,dy.$$
After summation, we deduce that:
$$e_N^\rho(\mu,\rho)\le \rho\int_{-\infty}^0(-y)^{\rho-1}\left(\frac{1}{N}\wedge F(y)\right)\,dy+\rho\int_0^{+\infty}y^{\rho-1}\left(\frac{1}{N}\wedge(1-F(y))\right)\,dy.$$
Since, by Fubini's theorem, $\displaystyle \rho\int_{-\infty}^0(-y)^{\rho-1}F(y)\,dy+\rho\int_0^{+\infty}y^{\rho-1}(1-F(y))\,dy=\int_\R |x|^\rho\mu(dx)<+\infty$, Lebesgue's theorem ensures that the right-hand side and therefore $e_N(\mu,\rho)$ go to $0$ as $N\to+\infty$.
\end{proof}

\section{The non compactly supported case}\label{parnoncomp}

According to Theorem 5.21 (ii) \cite{xuberger}, when the support of $\mu$ is bounded, $\displaystyle \sup_{N\ge 1}N^{1/\rho}e_N(\mu,\rho)<+\infty$ with $\displaystyle \lim_{N\to+\infty}N^{1/\rho}e_N(\mu,\rho)=0$ if and only if the quantile function $F^{-1}$ is continuous according to Remark 5.22 (ii) \cite{xuberger} and Proposition 2.1 \cite{BJ}. 
The case $\beta>1$ in the next example illustrates the possibility that, when $\rho>1$, $\displaystyle \lim_{N\to+\infty}N^{\frac{1}{\rho}}e_N(\mu,\rho)=0$ for some non compactly supported probability measures $\mu$. Of course, $F^{-1}$ is then continuous on $(0,1)$, since, by Remark 5.22 (ii) \cite{xuberger}, $\displaystyle \limsup_{N\to+\infty} N^{1/\rho}e_N(\mu,\rho)>0$ otherwise.
\begin{exple}\label{exempleexp}
For $\mu_\beta(dx)=f(x)\,dx$ with $f(x)=\mathbf{1}_{\{x>0\}}\beta x^{\beta -1}\exp\left(-x^{\beta}\right)$ with $\beta>0$ (the exponential distribution case $\beta=1$, was addressed in Example 5.17 and Remark 5.22 (i) \cite{xuberger}), we have that $F(x) = \mathbf{1}_{\{x>0\}}\left(1-\exp(-x^{\beta})\right)$, $F^{-1}(u) = \left(-\ln(1-u)\right)^{\frac{1}{\beta}}$ and $f\left(F^{-1}(u)\right)=\beta(1-u)(-\ln(1-u))^{1-\frac{1}{\beta}}$. The density $f$ is decreasing on $\left[x_\beta,+\infty\right)$ where $x_\beta=\left(\frac{(\beta-1)\vee 0}{\beta}\right)^{\frac{1}{\beta}}$. Using \eqref{minotermbordlimB}, the equality $F^{-1}(w)-F^{-1}(u)=\int_u^w\frac{dv}{f\left(F^{-1}(v)\right)}$ valid for $u,w\in(0,1)$ and the monotonicity of the density, we obtain that for $N$ large enough so that $\lceil NF(x_\beta)\rceil\le N-1$, 
\begin{align}
  e_N^\rho(\mu_\beta,\rho)&\ge \frac{1}{4^{\rho}N}\sum_{i=\lceil NF(x_\beta)\rceil +1}^N\left(\int_{\frac{4i-3}{4N}}^{\frac{4i-1}{4N}}\frac{du}{f(F^{-1}(u))}\right)^\rho\ge \frac{1}{8^{\rho}N^{\rho+1}}\sum_{i=\lceil NF(x_\beta)\rceil +1}^N\frac{1}{f^\rho\left(F^{-1}\left(\frac{4i-3}{4N}\right)\right)}\notag\\
  &\ge \frac{1}{(8N)^{\rho}}\sum_{i=\lceil NF(x_\beta)\rceil +2}^N\int_{\frac{i-2}{N}}^{\frac{i-1}{N}}\frac{du}{f^\rho\left(F^{-1}(u)\right)}=\frac{1}{(8N)^{\rho}}\int_{\frac{\lceil NF(x_\beta)\rceil}{N}}^{\frac{N-1}{N}}\frac{du}{f^\rho\left(F^{-1}(u)\right)}.\label{minoenrd}
\end{align}
Using H\"older's inequality for the second inequality, then Fubini's theorem for the third, we obtain that
\begin{align*}
  e_N^\rho(\mu_\beta,\rho)&-\int_{\frac{N-1}{N}}^{1}\left|x_N^N-F^{-1}(u)\right|^\rho\,du \le \sum_{i=1}^{N-1}\int_{\frac{i-1}{N}}^{\frac{i}{N}}\left|\int_{\frac{2i-1}{2N}}^u\frac{dv}{f(F^{-1}(v))}\right|^\rho\,du \\&\le \sum_{i=1}^{N-1}\int_{\frac{i-1}{N}}^{\frac{i}{N}}\left|u-\frac{2i-1}{2N}\right|^{\rho-1}\left|\int_{\frac{2i-1}{2N}}^u \frac{dv}{f^{\rho}\left(F^{-1}(v)\right)}\right|\,du\le\frac{1}{(2N)^\rho\rho}\int_0^{\frac{N-1}{N}}\frac{dv}{f^\rho\left(F^{-1}(v)\right)}.
\end{align*}
We have $F(x_\beta)<1$ and, when $\beta>1$, $F(x_\beta)>0$. By integration by parts, for $\rho>1$,
\begin{align*}
  (\rho-1)&\int_{F(x_\beta)}^{\frac{N-1}{N}}\frac{\beta^\rho\,du}{f^\rho\left(F^{-1}(u)\right)} = \int_{F(x_\beta)}^{\frac{N-1}{N}}(\rho-1)(1-u)^{-\rho}(-\ln(1-u))^{\frac{\rho}{\beta}-\rho}\,du\\
  &=\left[(1-u)^{1-\rho}(-\ln(1-u))^{\frac{\rho}{\beta}-\rho}\right]^{\frac{N-1}{N}}_{F(x_\beta)}+{\left(\frac{\rho}{\beta}-\rho\right)}\int_{F(x_\beta)}^{\frac{N-1}{N}}(1-u)^{-\rho}(-\ln(1-u))^{\frac{\rho}{\beta}-\rho-1}\,du\\
  &=N^{\rho-1}(\ln N)^{\frac{\rho}{\beta}-\rho}+o\left(\int_{F(x_\beta)}^{\frac{N-1}{N}}(1-u)^{-\rho}(-\ln(1-u))^{\frac{\rho}{\beta}-\rho}\,du\right)\sim N^{\rho-1}(\ln N)^{\frac{\rho}{\beta}-\rho},
\end{align*}
as $N\to+\infty$. We obtain the same equivalent when replacing the lower integration limit $F(x_\beta)$ in the left-hand side by $\frac{\lceil NF(x_\beta)\rceil}{N}$ or $0$ since $\displaystyle \lim_{N\to+\infty}\int^{\frac{\lceil NF(x_\beta)\rceil}{N}}_{F(x_\beta)}\frac{du}{f^\rho(F^{-1}(u))}=0$  and $\displaystyle \int_0^{F(x_\beta)}\frac{du}{f^\rho(F^{-1}(u))}<+\infty$.
On the other hand,
\begin{align*}
  \int_{\frac{N-1}{N}}^1\left|x_N^N-F^{-1}(u)\right|^\rho\,du\le \int_{\frac{N-1}{N}}^1\left(\left(-\ln(1-u)\right)^{\frac{1}{\beta}}-(\ln N)^{\frac{1}{\beta}}\right)^\rho\,du.
\end{align*}
When $\beta<1$, for $u\in\left[\frac{N-1}{N},1\right]$, $\left(-\ln(1-u)\right)^{\frac{1}{\beta}}-(\ln N)^{\frac{1}{\beta}} \le {\frac{1}{\beta}}\left(-\ln(1-u)\right)^{\frac{1}{\beta}-1}\left(-\ln(1-u)-\ln N\right)$ so that
\begin{align*}
   \int_{\frac{N-1}{N}}^1\left(\left(-\ln(1-u)\right)^{\frac{1}{\beta}}-(\ln N)^{\frac{1}{\beta}}\right)^\rho\,du&\le \frac{1}{\beta^\rho}\int_{\frac{N-1}{N}}^1\left(-\ln(1-u)\right)^{\frac{\rho}{\beta}-\rho}\left(-\ln(N(1-u)) \right)^\rho\,du\\
   &=\frac{1}{\beta^\rho N}\int_0^1\left(\ln N-\ln v\right)^{\frac{\rho}{\beta}-\rho}(-\ln(v))^\rho\,dv\\
   &\le \frac{2^{(\frac{\rho}{\beta}-\rho-1)\vee 0}}{\beta^\rho N}\left(\left(\ln N\right)^{\frac{\rho}{\beta}-\rho}\int_0^1(-\ln(v))^\rho\,dv+\int_0^1(-\ln(v))^{\frac{\rho}{\beta}}\,dv\right).
\end{align*}
When $\beta\ge 1$, for $N\ge 2$ and $u\in\left[\frac{N-1}{N},1\right]$, $\left(-\ln(1-u)\right)^{\frac{1}{\beta}}-(\ln N)^{\frac{1}{\beta}}\le {\frac{1}{\beta}}\left(\ln N\right)^{\frac{1}{\beta}-1}\left(-\ln(1-u)-\ln N\right)$ so that
\begin{align}
  \int_{\frac{N-1}{N}}^1\left(\left(-\ln(1-u)\right)^{\frac{1}{\beta}}-(\ln N)^{\frac{1}{\beta}}\right)^\rho\,du\le \frac{\left(\ln N\right)^{\frac{\rho}{\beta}-\rho}}{\beta^\rho N}\int_0^1(-\ln(v))^{\frac{\rho}{\beta}}\,dv.\label{amjotailexp}
\end{align}
We conclude that for $\rho>1$ and $\beta>0$, $e_N(\mu_\beta,\rho)\asymp N^{-\frac{1}{\rho}}(\ln N)^{\frac{1}{\beta}-1}\asymp {\cal W}_\rho(\frac{1}{N}(\sum_{i=1}^{N-1}\delta_{F^{-1}\left(\frac{2i-1}{2N}\right)}+\delta_{F^{-1}\left(\frac{N-1}{N}\right)}),\mu_\beta)$. In view of Theorem 5.20 \cite{xuberger}, this rate of convergence does not extend continuously to $e_N(\mu_\beta,1)$, at least for $\beta>1$. Indeed, by Remark 2.2 \cite{jrdcds}, $e_N(\mu_\beta,1)\asymp N^{-1}(\ln N)^{\frac{1}{\beta}}$, which in view of \eqref{amjotailexp}, implies that $\sum\limits_{i=1}^{N-1}\int_{\frac{i-1}{N}}^{\frac{i}{N}}\left|x_i^N-F^{-1}(u)\right|^\rho\,du\asymp N^{-1}(\ln N)^{\frac{1}{\beta}}$. 

In the Gaussian tail case $\beta=2$, $e_N(\mu_2,\rho)\asymp N^{-\frac{1}{\rho}}(\ln N)^{-\frac{1}{2}+{\mathbf 1}_{\{\rho=1\}}}$ for $\rho\ge 1$, like for the true Gaussian distribution, according to Example 5.18 \cite{xuberger}. This matches the rate obtained when $\rho>2$ in Corollary 6.14 \cite{bobkovledoux} for $\E^{1/\rho}\left[{\cal W}_\rho^\rho\left(\frac{1}{N}\sum\limits_{i=1}^N\delta_{X_i},\mu\right)\right]$ where $(X_i)_{i\ge 1}$ are i.i.d. with respect to some Gaussian distribution $\mu$ with positive variance. When $\rho=2$, still according to this corollary the random rate is $N^{-1/2}(\ln\ln N)^{1/2}$ (of course worse than the standard Monte Carlo rate $N^{-1/2}$).\end{exple}
According to the next result, the order of convergence of $e_N(\mu,\rho)$ cannot exceed $\frac{1}{\rho}$ when the support of $\mu$ is not bounded.
\begin{prop}
 Let $\rho>1$. Then $\exists \, \alpha>\frac 1\rho,\;\sup_{N\ge 1}N^\alpha e_N(\mu,\rho)<+\infty\implies F^{-1}(1)-F^{-1}(0+)<+\infty$. \label{propals1rcomp}
\end{prop}
\begin{proof}
Let $\rho>1$ and $\alpha>\frac 1\rho$ be such that $\sup_{N\ge 1}N^{\alpha} e_N(\mu,\rho)<+\infty$ so that, by \eqref{enrho}, 
$$\sup_{N\ge 1}N^{\alpha\rho}\left(\int_0^{\frac 1 N} \left|F^{-1}(u)-x^N_1\right|^{\rho}\,du+\int^1_{\frac{N-1}{N}} \left|F^{-1}(u)-x^N_N\right|^{\rho}\,du\right)<+\infty.$$
By \eqref{minotermbordlimB} for $i=1$ and $N\ge1$, we have: 
\begin{align*}
 \displaystyle \int_0^{\frac 1 N} \left|F^{-1}(u)-x^N_1\right|^{\rho}\,du \ge \frac{1}{4^\rho N}\left(F^{-1}\left(\frac{1}{2N}\right)-F^{-1}\left(\frac{1}{4N}\right)\right)^{\rho}.
\end{align*}
Therefore  $C:=\sup\limits_{N\ge 1}(2N)^{\alpha-\frac{1}{\rho}}\left(F^{-1}\left(\frac{1}{2N}\right)-F^{-1}\left(\frac{1}{4N}\right)\right)<+\infty$. For $k\in \N^{*}$, we deduce that $F^{-1}\left(2^{-(k+1)}\right)-F^{-1}\left(2^{-k}\right) \ge - C2^{\frac{1-\alpha \rho}{\rho}k}$, and after summation that:
\begin{equation}
   \forall k\in \N^{*},F^{-1}\left(2^{-k}\right) \ge F^{-1}\left(1/2\right) - \frac{C}{2^{\alpha -\frac{1}{\rho}} -1 }\left(1- 2^{\frac{1-\alpha \rho}{\rho}(k-1)}\right).\label{minoF-1u0}
\end{equation}
When $k\to+\infty$, the right-hand side goes to $\left(F^{-1}\left(\frac 1 2\right) - \frac{C}{2^{\alpha -\frac{1}{\rho}}-1}\right)>-\infty$ so that $F^{-1}(0+)>-\infty$. In a symmetric way, we check that $F^{-1}(1)<+\infty$ so that $\mu$ is compactly supported.
\end{proof}
The next theorem gives a necessary and sufficient condition for $e_N(\mu,\rho)$ to go to $0$ with order $\alpha\in\left(0,\frac 1\rho\right)$.
\begin{thm}\label{alphaRater}
Let $\rho\ge 1$ and $\alpha\in \left(0,\frac{1}{\rho}\right)$. We have \begin{align*}
  {\sup_{x\ge 0}}\;x^{\frac{\rho}{1-\alpha \rho}}\Big(F(-x)+1-F(x)\Big)<+\infty \Leftrightarrow \displaystyle{\sup_{N \ge 1}} N^{\alpha} \, e_N(\mu,\rho)<+\infty\Leftrightarrow \displaystyle{\sup_{N \ge 2}}\sup_{x_{2:N-1}}N^{\alpha}{\cal W}_\rho(\mu_N(x_{2:N-1}),\mu)<+\infty\end{align*}
where $\mu_N(x_{2:N-1})=\frac{1}{N}\left(\delta_{F^{-1}\left(\frac{1}{N}\right)\wedge (-N^{\frac{1}{\rho}-\alpha})}+\sum_{i=2}^{N-1}\delta_{x_i}+\delta_{F^{-1}\left(\frac{N-1}{N}\right)\vee N^{\frac{1}{\rho}-\alpha}}\right)$ and $\sup_{x_{2:N-1}}$ means the supremum over the choice of $x_i\in[F^{-1}(\frac{i-1}{N}+),F^{-1}(\frac{i}{N})]$ for $2\le i\le N-1$. Moreover,
$$\lim_{x\to +\infty}x^{\frac{\rho}{1-\alpha \rho}}\Big(F(-x)+1-F(x)\Big)=0\Leftrightarrow \lim_{N\to+\infty}N^{\alpha}e_N(\mu,\rho)=0.$$

\end{thm}

Let $\hat\beta=\sup\{\beta\ge 0:\int_\R|x|^\beta\mu(dx)<\infty\}$. For $\beta \in (0,\hat\beta)$, $\int_\R|x|^\beta\mu(dx)<\infty$ while when $\hat\beta<\infty$, for $\beta>\hat\beta$, $\int_\R|x|^{\frac{\hat\beta+\beta}{2}}=+\infty$, so that by Lemma \ref{lemququcdf} ${\sup_{x\ge 0}}\;x^{\beta}\Big(F(-x)+1-F(x)\Big)=+\infty$. Let $\rho\ge 1$. If $\hat\beta>\rho$, we deduce from Theorem 5.21 (i) \cite{xuberger} that for each $\alpha\in (0,\frac 1\rho-\frac 1 {\hat\beta})$, $\lim_{N\to\infty}N^\alpha e_N(\mu,\rho)=0$ and, when $\hat\beta<+\infty$, Theorem \ref{alphaRater} ensures that for each $\alpha>\frac 1\rho-\frac 1 {\hat\beta}$, $\sup_{N\ge 1} N^\alpha e_N(\mu,\rho)=+\infty$ since $\frac\rho{1-\alpha\rho}>\hat\beta$. In this sense, when $\rho<\hat\beta<+\infty$ the order of convergence of $e_N(\mu,\rho)$ to $0$ is $\frac 1\rho-\frac 1 {\hat\beta}$. Moreover, the boundedness and the vanishing limit at infinity for the sequence $(N^{\frac 1\rho-\frac 1 {\hat\beta}} e_N(\mu,\rho))_{N\ge 1}$ are respectively equivalent to the same property for the function $\R_+\ni x\mapsto x^{\hat\beta}\Big(F(-x)+1-F(x)\Big)<+\infty$. Note that $\limsup_{x\to+\infty} x^{\hat\beta}\Big(F(-x)+1-F(x)\Big)$ can be $0$, or positive or $+\infty$.

\begin{remark}\begin{itemize}
 \item In the proof we check that if $C={\sup\limits_{x\ge 0}}\;x^{\frac{\rho}{1-\alpha \rho}}\Big(F(-x)+1-F(x)\Big)$, then \begin{align*}
  \sup_{N\ge 1} N^{\alpha\rho}e^\rho_N(\mu,\rho)\le \sup_{N \ge 2}\sup_{x_{2:N-1}}N^{\alpha\rho}{\cal W}^\rho_\rho(\mu_N(x_{2:N-1}),\mu)\le 2^\rho C^{1-\alpha\rho}+\frac{1-\alpha\rho}{\alpha\rho}C+1+2^{\rho-1}+2^{\rho+\alpha\rho-2}\left|F^{-1}\left(1/2\right)\right|^\rho.
\end{align*}
 \item According to Theorem 7.16 \cite{bobkovledoux}, for $(X_i)_{i\ge 1}$ i.i.d. according to $\mu$, $$\sup_{N\ge 1}N^{\frac{1}{2\rho}}\E^{1/\rho}\left[{\cal W}_\rho^\rho\left(\frac{1}{N}\sum_{i=1}^N\delta_{X_i},\mu\right)\right]\le \left(\rho 2^{\rho-1}\int_\R |x|^{\rho-1}\sqrt{F(x)(1-F(x))}\,dx\right)^{1/\rho}$$ with $\exists \;\varepsilon>0,\;\int_\R|x|^{2\rho+\varepsilon}\mu(dx)<+\infty\Rightarrow \int_\R |x|^{\rho-1}\sqrt{F(x)(1-F(x))}\,dx<+\infty\Rightarrow \int_\R|x|^{2\rho}\mu(dx)<+\infty$ by the discussion just after the theorem. The condition ${\sup_{x\ge 0}}\;x^{2\rho}\Big(F(-x)+1-F(x)\Big)<+\infty$ equivalent to $\sup_{N \ge 1} N^{\frac{1}{2\rho}} \, e_N(\mu,\rho)<+\infty$ is slightly weaker, according to Lemma \ref{lemququcdf} just below. Moreover, we address similarly any order of convergence $\alpha$ with $\alpha\in \left(0,\frac{1}{\rho}\right)$ for $e_N(\mu,\rho)$, while the order $\frac{1}{2\rho}$ seems to play a special role for $\E^{1/\rho}\left[{\cal W}_\rho^\rho\left(\frac{1}{N}\sum\limits_{i=1}^N\delta_{X_i},\mu\right)\right]$ in the random case. When $\rho=1$, the order of convergence $\alpha$ for $\alpha\in(0,1/2)$ is addressed in the random case in Theorem 2.2 \cite{barrGinMat} where the finiteness of ${\sup_{x\ge 0}}\,x^{\frac{1}{1-\alpha}}\Big(F(-x)+1-F(x)\Big)$ is stated to be equivalent to the stochastic boundedness of the sequence $\left(N^\alpha{\cal W}_1\left(\frac{1}{N}\sum\limits_{i=1}^N\delta_{X_i},\mu\right)\right)_{N\ge 1}$. When $\alpha=1/2$, the stochastic boundedness property is, according to Theorem 2.1 (b) \cite{barrGinMat}, equivalent to $\int_\R\sqrt{F(x)(1-F(x))}\,dx<+\infty$.
\end{itemize}
\end{remark}  
The proof of Theorem \ref{alphaRater} relies on the two next lemmas.
\begin{lem}\label{lemququcdf}
For $\beta>0$, we have 
   \begin{align*}
      \int_{\R}|y|^{\beta}\mu(dy)<+\infty &\implies \lim_{x\to+\infty}x^{\beta}\Big(F(-x)+1-F(x)\Big)=0 \\&\implies {\sup_{x\ge 0}}\;x^{\beta}\Big(F(-x)+1-F(x)\Big)<+\infty \implies \forall \varepsilon\in \left(0,\beta\right],\;\int_{\R}|y|^{\beta-\varepsilon}\mu(dy)<+\infty
   \end{align*}
   and ${\sup_{x\ge 0}}\;x^{\beta}\Big(F(-x)+1-F(x)\Big)<+\infty \Leftrightarrow \sup_{u\in(0,1/2]}u^{\frac{1}{\beta}}\left(F^{-1}(1-u)-F^{-1}(u)\right)<+\infty$ with
   \begin{equation}
   \sup_{u\in(0,1/2]}u^{\frac{1}{\beta}}\left(F^{-1}(1-u)-F^{-1}(u)\right) \le \left(\sup_{x\ge 0}\,x^{\beta}F(-x)\right)^{\frac{1}{\beta}}+\left(\sup_{x\ge 0}\,x^{\beta}(1-F(x))\right)^{\frac{1}{\beta}}.\label{majodifquantb}
 \end{equation}
 Last, $\lim_{x\to +\infty}x^{\beta}\Big(F(-x)+1-F(x)\Big)=0\Leftrightarrow \lim_{u\to 0+}u^{\frac{1}{\beta}}\left(F^{-1}(1-u)-F^{-1}(u)\right)=0$.
\end{lem}
\begin{proof}
  Let $\beta>0$. For $x>0$, using the monotonicity of $F$ for the first inequality then that for $y\in[\frac x 2,x]$, $y^{\beta-1}\ge \left(\frac{x}{2}\right)^{\beta-1}\wedge x^{\beta-1}=\frac{x^{\beta-1}}{2^{(\beta-1)\vee 0}}$, we obtain that 
\begin{align*}
  F(-x)+1-F(x) \le \frac{2}{x}\int_{x/2}^{x}\Big(F(-y)+1-F(y)\Big)\,dy &\le \frac{2^{\beta\vee 1}}{x^\beta}\int_{x/2}^{+\infty}y^{\beta-1}\Big(F(-y)+1-F(y)\Big)\,dy.\end{align*}
Since $\int_0^{+\infty}y^{\beta-1}\Big(F(-y)+1-F(y)\Big)\,dy=\frac 1\beta\int_{\R}|y|^{\beta}\mu(dy)$, the finiteness of $\int_{\R}|y|^{\beta}\mu(dy)$ implies by Lebesgue's theorem that $\lim_{x\to\infty} x^\beta\left(F(-x)+1-F(x)\right)=0$. Since $x\mapsto x^\beta\left(F(-x)+1-F(x)\right)$ is right-continuous with left-hand limits on $[0,+\infty)$, $${\sup_{x\ge 0}}\;x^{\beta}\Big(F(-x)+1-F(x)\Big)<+\infty \Leftrightarrow \limsup_{x\to\infty}x^{\beta}\Big(F(-x)+1-F(x)\Big)<+\infty,$$
with the latter property clearly implied by $\lim_{x\to\infty} x^\beta\left(F(-x)+1-F(x)\right)=0$.

For  $\varepsilon\in (0,\beta)$, using that for $y\ge 0$, $F(-y)+1-F(y)=\mu((-\infty,-y]\cup (y,+\infty))\le 1$, we obtain
  \begin{align*}
   \int_\R|x|^{\beta-\varepsilon}\mu(dx)&=(\beta-\varepsilon)\int_0^{+\infty}y^{\beta-\varepsilon-1}(F(-y)+1-F(y))\,dy\\
   &\le (\beta-\varepsilon)\int_0^1y^{\beta-\varepsilon-1}dy+(\beta-\varepsilon)\sup\limits_{x \ge 0}x^{\beta}\Big(F(-x)+1-F(x)\Big)\int_1^{+\infty}y^{-\varepsilon-1}\,dy\\
   &=1+\frac{\beta-\varepsilon}{\varepsilon}\sup\limits_{x \ge 0}x^{\beta}\Big(F(-x)+1-F(x)\Big).
  \end{align*}
Therefore $\displaystyle \sup\limits_{x \ge 0}x^{\beta}\Big(F(-x)+1-F(x)\Big)<+\infty\implies \forall \varepsilon\in (0,\beta),\;\int_\R|x|^{\beta-\varepsilon}\mu(dx)<+\infty$.\\
Let us next check that \begin{equation}
   \sup_{x\ge 0}\,x^{\beta}\Big(F(-x)+(1-F(x))\Big)<+\infty\Leftrightarrow \sup_{u\in(0,1/2]}u^{\frac{1}{\beta}}(F^{-1}(1-u)-F^{-1}(u))<+\infty\label{equivtailquant}.
 \end{equation} For the necessary condition, we set $u\in(0,1/2]$. Either $F^{-1}(u)\ge 0$ or, since for all $v\in (0,1)$, $F(F^{-1}(v))\ge v$, we have $\left(-F^{-1}(u)\right)^\beta u\le \sup_{x\ge -F^{-1}(u)}\,x^{\beta}F(-x)$ and therefore $F^{-1}(u)\ge -\left(\sup_{x\ge 0}\,x^{\beta}F(-x)\right)^{\frac{1}{\beta}}u^{-\frac{1}{\beta}}$. Either $F^{-1}(1-u)\le 0$ or, since for all $v\in (0,1)$, $F(F^{-1}(v)-)\le v$,  we have $(F^{-1}(1-u))^\beta u\le \sup_{x\ge F^{-1}(1-u)}\,x^{\beta}(1-F(x-))$ and therefore $F^{-1}(1-u)\le \left(\sup_{x\ge 0}\,x^{\beta}(1-F(x))\right)^{\frac{1}{\beta}}u^{-\frac{1}{\beta}}$. Hence \eqref{majodifquantb} holds.
 
 For the sufficient condition, we remark that the finiteness of $\sup_{u\in(0,1/2]}u^{\frac{1}{\beta}}(F^{-1}(1-u)-F^{-1}(u))$ and the inequality $F^{-1}(1-u)-F^{-1}(u)\ge \left(F^{-1}(1/2)-F^{-1}(u)\right)\vee\left(F^{-1}(1-u)-F^{-1}(1/2)\right) $ valid for $u\in(0,1/2]$ imply that $\inf_{u \in (0,1/2]}u^{\frac{1}{\beta}}F^{-1}(u) >-\infty$ and $\sup_{u \in (0,1/2]}u^{\frac{1}{\beta}}F^{-1}(1-u)<+\infty$. With the inequality $ x\ge F^{-1}(F(x))$ valid for $x\in\R$ such that $0<F(x)<1$, this implies that
 $\inf_{x\in\R:F(x)\le 1/2}\left(F(x)\right)^{\frac{1}{\beta}}x> -\infty$ and therefore that $\sup_{x\ge 0}x^{{\beta}}F(-x)<+\infty$. With the inequality $ x\le F^{-1}(F(x)+)$ valid for $x\in\R$ such that $0<F(x)<1$, we obtain, in a symmetric way $\sup_{x\ge 0}x^{{\beta}}(1-F(x))<+\infty$.

 Let us finally check that $\lim_{x\to +\infty}x^{\beta}\Big(F(-x)+1-F(x)\Big)=0\Leftrightarrow \lim_{u\to 0+}u^{\frac{1}{\beta}}\left(F^{-1}(1-u)-F^{-1}(u)\right)=0$. For the necessary condition, we remark that either $F^{-1}(1)<+\infty$ and $\lim_{u\to 0+}u^{\frac 1\beta}F^{-1}(1-u)=0$ or $F^{-1}(1-u)$ goes to $+\infty$ as $u\to 0+$. For $u$ small enough so that $F^{-1}(1-u)>0$, we have $(F^{-1}(1-u))^\beta u\le \sup_{x\ge F^{-1}(1-u)}\,x^{\beta}(1-F(x-))$, from which we deduce that $\lim_{u\to 0+}u(F^{-1}(1-u))^\beta =0$. The fact that $\lim_{u\to 0+}u^{\frac 1\beta}F^{-1}(u)=0$ is deduced by a symmetric reasoning.

 For the sufficient condition, we use that $x(1-F(x))^{\frac{1}{\beta}}\le \sup_{u\le 1-F(x)}u^{\frac{1}{\beta}}F^{-1}((1-u)+)$ and $x\left(F(x)\right)^{\frac{1}{\beta}}\ge \inf_{u\le F(x)}u^{\frac{1}{\beta}}F^{-1}(u)$ for $x\in\R$ such that $0<F(x)<1$.

\end{proof}

\begin{lem}\label{lemcontxn}
  Let $\rho\ge 1$ and $\alpha\in \left(0,\frac{1}{\rho}\right)$.There is a finite constant $C$ only depending on $\rho$ and $\alpha$ such that the two extremal points in the optimal sequence $(x_i^N)_{1\le i\le N}$ for $e_N(\mu,\rho)$ satisfy
  $$\forall N\ge 1,\;x_1^N\ge CN^{\frac{1}{\rho}-\alpha}\inf_{u\in (0,\frac 1 N)}u^{\frac{1}{\rho}-\alpha}F^{-1}(u)\mbox{ and }x_N^N\le CN^{\frac{1}{\rho}-\alpha}\sup_{u\in (0,\frac 1 N)}u^{\frac{1}{\rho}-\alpha}F^{-1}(1-u).$$

  If $\sup_{u\in(0,1/2]}u^{\frac{1}{\rho}-\alpha}\left(F^{-1}(1-u)-F^{-1}(u)\right)<+\infty$, then $\sup_{N\ge 1}N^{\alpha-\frac{1}{\rho}}\left(x_N^N\vee \left(-x_1^N\right)\right)<+\infty$.
\end{lem}

\begin{proof}
  Since the finiteness of $\sup_{u\in(0,1/2]}u^{\frac{1}{\rho}-\alpha}\left(F^{-1}(1-u)-F^{-1}(u)\right)$ implies the finiteness of both $\sup_{u\in (0,1)}u^{\frac{1}{\rho}-\alpha}F^{-1}(1-u)$ and $\inf_{u\in (0,1)}u^{\frac{1}{\rho}-\alpha}F^{-1}(u)$, the second statement is a consequence of the first one, that we are now going to prove.
When $\rho=1$ (resp. $\rho=2$), then the conclusion easily follows from the explicit form $x_1^N=F^{-1}\left(\frac{1}{2N}\right)$ and $x_N^N=F^{-1}\left(\frac{2N-1}{2N}\right)$ (resp. $x_1^N=N\int_0^{\frac{1}{N}}F^{-1}(u)\,du$ and $x_N^N=N\int^1_{\frac{N-1}{N}}F^{-1}(u)\,du$). In the general case $\rho>1$, we are going to take advantage of the expression 
$$f(y)=\rho\int_0^{\frac{1}{N}}\left({\mathbf 1}_{\{y\ge F^{-1}(1-u)\}}\left(y-F^{-1}(1-u)\right)^{\rho-1}-{\mathbf 1}_{\{y< F^{-1}(1-u)\}}\left(F^{-1}(1-u)-y\right)^{\rho-1}\right)\,du$$
of the derivative of the function $\displaystyle \R\ni y\mapsto \int_{\frac{N-1}{N}}^{1}\left|y-F^{-1}(u)\right|^\rho\,du$ minimized by $x_N^N$. Since this function is strictly convex $x_N^N=\inf\{y\in\R: f(y)\ge 0\}$.
Let us first suppose that $S_N:=\sup_{u\in (0,\frac 1 N)}u^{\frac{1}{\rho}-\alpha}F^{-1}(1-u)\in (0,+\infty)$. Since for fixed $y\in\R$,
$\R\ni x\mapsto \left({\mathbf 1}_{\{y\ge x\}}(y-x)^{\rho-1}-{\mathbf 1}_{\{y< x\}}(x-y)^{\rho-1}\right)$ is non-increasing, we deduce that $\forall y\in\R$, $f(y)\ge \rho S_N^{\rho-1}g(\frac{y}{S_N})$ where
$$g(z)=\int_0^{\frac{1}{N}}\left({\mathbf 1}_{\left\{z\ge u^{\alpha-\frac{1}{\rho}}\right\}}\left(z-u^{\alpha-\frac{1}{\rho}}\right)^{\rho-1}-{\mathbf 1}_{\left\{z< u^{\alpha-\frac{1}{\rho}}\right\}}\left(u^{\alpha-\frac{1}{\rho}}-z\right)^{\rho-1}\right)\,du.$$
For $z\ge (4N)^{\frac{1}{\rho}-\alpha}$, we have $z^{\frac{\rho}{\alpha\rho-1}}\le\frac{1}{4N}$ and $z-(2N)^{\frac{1}{\rho}-\alpha}\ge \left(1-2^{\alpha-\frac 1\rho}\right)z$ so that 
\begin{align*}
  g(z)&=\int_{z^{\frac{\rho}{\alpha\rho-1}}}^{\frac{1}{N}}\left(z-u^{\alpha-\frac{1}{\rho}}\right)^{\rho-1}\,du-\int_0^{z^{\frac{\rho}{\alpha\rho-1}}}\left(u^{\alpha-\frac{1}{\rho}}-z\right)^{\rho-1}\,du \ge \int_{\frac{1}{2N}}^{\frac{1}{N}}\left(z-(2N)^{\frac{1}{\rho}-\alpha}\right)^{\rho-1}\,du-\int_0^{z^{\frac{\rho}{\alpha\rho-1}}}u^{(\rho-1)\frac{\alpha\rho-1}{\rho}}\,du\\
  &\ge \left(1-2^{\alpha-\frac 1\rho}\right)^{\rho-1}z^{\rho-1}\int_{\frac{1}{2N}}^{\frac{1}{N}}\,du-\frac{\rho z^{\frac{\rho}{\alpha\rho-1}+\rho-1}}{1+(\rho-1)\alpha\rho} = z^{\rho-1}\left(\frac{\left(1-2^{\alpha-\frac 1\rho}\right)^{\rho-1}}{2N}-\frac{\rho z^{\frac{\rho}{\alpha\rho-1}}}{1+(\rho-1)\alpha\rho} \right).
\end{align*}
The right-hand side is positive for $z>(\kappa N)^{\frac{1}{\rho}-\alpha}$ with $\kappa:=\frac{2\rho}{\left(1-2^{\alpha-\frac 1\rho}\right)^{\rho-1}\left(1+(\rho-1)\alpha\rho\right)}$. Hence for $z>\left(\left(\kappa\vee 4\right)N\right)^{\frac{1}{\rho}-\alpha}$, $g(z)>0$ so that for $y>\left(\left(\kappa\vee 4\right)N\right)^{\frac{1}{\rho}-\alpha}S_N$, $f(y)>0$ and therefore $$x_N^N \le \left(\left(\kappa\vee 4\right)N\right)^{\frac{1}{\rho}-\alpha}S_N .$$
Clearly, this inequality remains valid when $S_N=+\infty$. It holds in full generality since  $S_N\ge 0$ and $S_N=0\Leftrightarrow F^{-1}(1)\le 0$ a condition under which $x_N^N\le 0$ since $x_N^N\in\left[F^{-1}\left(\frac{N-1}{N}+\right),F^{-1}\left(1\right)\right]\cap\R$. 
By a symmetric reasoning, we check that $x_1^N\ge \left(\left(\kappa\vee 4\right)N\right)^{\frac{1}{\rho}-\alpha}\inf_{u\in (0,\frac 1 N)}u^{\frac{1}{\rho}-\alpha}F^{-1}(u)$.
\end{proof}

\begin{proof}[Proof of Theorem \ref{alphaRater}]
Since by Lemma \ref{lemququcdf}, $$\sup_{u\in(0,1/2]}u^{\frac{1}{\rho}-\alpha}\left(F^{-1}(1-u)-F^{-1}(u)\right)<+\infty\implies {\sup_{x\ge 0}}\;x^{\frac{\rho}{1-\alpha \rho}}\Big(F(-x)+1-F(x)\Big)<+\infty$$ and, by \eqref{enrho}, $$e_N^\rho(\mu,\rho)\ge \int_0^{\frac 1 N} \left|F^{-1}(u)-x^N_1\right|^{\rho}\,du+\int^1_{\frac{N-1}{N}} \left|F^{-1}(u)-x^N_N\right|^{\rho}\,du$$ to prove the equivalence, it is enough to check that 
\begin{align*}
  &{\sup_{x\ge 0}}\;x^{\frac{\rho}{1-\alpha \rho}}\Big(F(-x)+1-F(x)\Big)<+\infty \implies \displaystyle{\sup_{N \ge 1}} N^{\alpha\rho} \, e_N^\rho(\mu,\rho)<+\infty\mbox{ and that }\\
  &\sup_{N\ge 1}N^{\alpha\rho}\left(\int_0^{\frac 1 N} \left|F^{-1}(u)-x^N_1\right|^{\rho}\,du+\int^1_{\frac{N-1}{N}} \left|F^{-1}(u)-x^N_N\right|^{\rho}\,du\right)<+\infty\\&\phantom{\sup_{N\ge 1}N^{\alpha\rho}\int_0^{\frac 1 N} \left|F^{-1}(u)-x^N_1\right|^{\rho}\,du}
\implies\sup_{u\in(0,1/2]}u^{\frac{1}{\rho}-\alpha}\left(F^{-1}(1-u)-F^{-1}(u)\right)<+\infty.
\end{align*}
We are now going to do so and thus prove that the four suprema in the two last implications are simultaneously finite or infinite.

Let us first suppose that $C:={\sup_{x\ge 0}}\;x^{\frac{\rho}{1-\alpha \rho}}\Big(F(-x)+1-F(x)\Big)<+\infty$ and set $N\ge 2$.
Let $x_i\in[F^{-1}(\frac{i-1}{N}+),F^{-1}(\frac{i}{N})]$ for $2\le i\le N-1$ 
. We have
$$e^\rho_N(\mu,\rho)\le {\cal W}^\rho_\rho\left(\mu_N(x_{2:N-1}),\mu\right)= L_N+M_N+U_N$$
with $L_N=\int_{0}^{\frac{1}{N}}\left|F^{-1}(u)-F^{-1}\left(\frac{1}{N}\right)\wedge (-N^{\frac{1}{\rho}-\alpha})\right|^\rho\,du$, $U_N=\int_{\frac{N-1}{N}}^{1}\left|F^{-1}(u)-F^{-1}\left(\frac{N-1}{N}\right)\vee N^{\frac{1}{\rho}-\alpha}\right|^\rho\,du$ and
\begin{align}
   M_N&=\sum_{i=2}^{N-1}\int_{\frac{i-1}{N}}^{\frac{i}{N}}\left|F^{-1}(u)-x_i\right|^\rho\,du \le \sum_{i=2}^{N-1}\int_{\frac{i-1}{N}}^{\frac{i}{N}}\left(F^{-1}\left(\frac{i}{N}\right)-F^{-1}\left(\frac{i-1}{N}\right)\right)^\rho\,du\notag\\&\le \frac{1}{N}\sum_{i=1}^N\left(F^{-1}\left(\frac{N-1}{N}\right)-F^{-1}\left(\frac{1}{N}\right)\right)^{\rho-1}\left(F^{-1}\left(\frac{i}{N}\right)-F^{-1}\left(\frac{i-1}{N}\right)\right)\notag\\&=\frac{1}{N}\left(F^{-1}\left(\frac{N-1}{N}\right)-F^{-1}\left(\frac{1}{N}\right)\right)^{\rho}\label{majomn}\\&\le 2^\rho C^{1-\alpha\rho}N^{-\alpha\rho},\notag\end{align}
where we used \eqref{majodifquantb} applied with $\beta=\frac{\rho}{1-\alpha \rho}$ for the last inequality. 
Let $x_+=0\vee x$ denote the positive part of any real number $x$. Applying Lemma \ref{lemenf} with $x=F^{-1}\left(\frac{1}{N}\right)\wedge \left(-N^{\frac{1}{\rho}-\alpha}\right)$, we obtain that
\begin{align*}
  L_N &= \rho\int_{-\infty}^{F^{-1}\left(\frac{1}{N}\right)\wedge \left(-N^{\frac{1}{\rho}-\alpha}\right)}\left(F^{-1}\left(\frac{1}{N}\right)\wedge \left(-N^{\frac{1}{\rho}-\alpha}\right)-y\right)^{\rho-1}F(y)\,dy\\&+\rho\int^{F^{-1}\left(\frac{1}{N}\right)}_{F^{-1}\left(\frac{1}{N}\right)\wedge \left(-N^{\frac{1}{\rho}-\alpha}\right)}\left(y-F^{-1}\left(\frac{1}{N}\right)\wedge \left(-N^{\frac{1}{\rho}-\alpha}\right)\right)^{\rho-1}\left(\frac{1}{N}-F(y)\right)\,dy
  \\&\le \rho\int^{+\infty}_{N^{\frac{1}{\rho}-\alpha}}y^{\rho-1}F(-y)\,dy+\frac{1}{N}\left(N^{\frac{1}{\rho}-\alpha}+F^{-1}\left(\frac{1}{N}\right)\right)_+^{\rho}.
\end{align*}
In a symmetric way, we check that $\displaystyle U_N\le \rho\int^{+\infty}_{N^{\frac{1}{\rho}-\alpha}}y^{\rho-1}(1-F(y))\,dy+\frac{1}{N}\left(N^{\frac{1}{\rho}-\alpha}-F^{-1}\left(\frac{N-1}{N}\right)\right)_+^{\rho}$ so that
\begin{align*}
   L_N+U_N &\le \rho C\int^{+\infty}_{N^{\frac{1}{\rho}-\alpha}}y^{-1-\frac{\alpha\rho^2}{1-\alpha\rho}}\,dy+\frac{1}{N}\left(\left(N^{\frac{1}{\rho}-\alpha}+F^{-1}\left(1/2\right)\right)_+^{\rho}+\left(N^{\frac{1}{\rho}-\alpha}-F^{-1}\left(1/2\right)\right)_+^{\rho}\right)\\
   &\le \frac{1-\alpha\rho}{\alpha\rho} CN^{-\alpha\rho}+\left(1+2^{\rho-1}\right)N^{-\alpha\rho}+2^{\rho-1}\left|F^{-1}\left(1/2\right)\right|^\rho N^{-1}.
\end{align*}
Since $N^{-1}\le 2^{\alpha\rho-1}N^{-\alpha\rho}$, we conclude that with $\sup_{x_{2:N-1}}$ denoting the supremum over $x_i\in[F^{-1}(\frac{i-1}{N}+),F^{-1}(\frac{i}{N})]$ for $2\le i\le N-1$,
\begin{align*}
   \sup_{N\ge 2}N^{\alpha\rho}e^\rho_N(\mu,\rho)\le \sup_{N\ge 2}\sup_{x_{2:N-1}}N^{\alpha\rho}{\cal W}_\rho^\rho(\mu_N(x_{2:N}),\mu)\le 2^\rho C^{1-\alpha\rho}+\frac{1-\alpha\rho}{\alpha\rho}C+1+2^{\rho-1}+2^{\rho+\alpha\rho-2}\left|F^{-1}\left(1/2\right)\right|^\rho.
\end{align*}
We may replace $\sup_{N\ge 2}N^{\alpha\rho}e^\rho_N(\mu,\rho)$ by $\sup_{N\ge 1}N^{\alpha\rho}e^\rho_N(\mu,\rho)$ in the left-hand side, since, applying Lemma \ref{lemenf} with $x=0$, then using that for $y\ge 0$, $F(-y)+1-F(y)=\mu((-\infty,-y]\cup (y,+\infty))\le 1$, we obtain that
\begin{align*}
   e^\rho_1(\mu,\rho) &\le \rho\int_0^{+\infty}y^{\rho-1}\Big(F(-y)+1-F(y)\Big)\,dy \le \rho\int_0^1y^{\rho-1}\,dy+\rho C\int^{+\infty}_{1}y^{-1-\frac{\alpha\rho^2}{1-\alpha\rho}}\,dy=1+\frac{1-\alpha\rho}{\alpha\rho}C.
\end{align*}
Let us next suppose that $\displaystyle \sup_{N\ge 1}N^{\alpha\rho}\left(\int_0^{\frac 1 N} \left|F^{-1}(u)-x^N_1\right|^{\rho}\,du+\int^1_{\frac{N-1}{N}} \left|F^{-1}(u)-x^N_N\right|^{\rho}\,du\right)<+\infty$. Like in the proof of Proposition \ref{propals1rcomp}, we deduce \eqref{minoF-1u0}. With the monotonicity of $F^{-1}$, this inequality implies that
$$\exists C<+\infty,\;\forall u \in (0,1/2], \quad F^{-1}(u) \ge F^{-1}(1/2) - \frac{C}{1 - 2^{\alpha-\frac{1}{\rho}}}\left(u^{\alpha-\frac{1}{\rho}}-1\right), $$
and therefore that $\inf_{u \in (0,1/2]}\left(u^{\frac{1}{\rho}-\alpha}F^{-1}(u)\right) >-\infty$. With a symmetric reasoning, we conclude that $$\sup_{u\in(0,1/2]}u^{\frac{1}{\rho}-\alpha}\Big(F^{-1}(1-u)-F^{-1}(u)\Big)<+\infty.$$

Let us now assume that $\limsup_{x\to+\infty}x^{\frac{\rho}{1-\alpha\rho}}\Big(F(-x)+1-F(x)\Big)\in(0,+\infty)$, which, in particular implies that $\sup_{x\ge 0}x^{\frac{\rho}{1-\alpha\rho}}\Big(F(-x)+1-F(x)\Big)<+\infty$ and check that $\limsup_{N\to\infty}N^\alpha e_N(\mu,\rho)>0$. For $x>0$, we have, on the one hand
\begin{align*}
   x^{\frac{\alpha \rho^2}{1-\alpha\rho}}\int_x^{+\infty}y^{\rho-1}\Big(F(-y)+1-F(y)\Big)\,dy &\ge x^{\frac{\alpha \rho^2}{1-\alpha\rho}}\int_x^{2x}x^{\rho-1}\Big(F(-2x)+1-F(2x)\Big)\,dy\\
   &=x^{\frac{\rho}{1-\alpha\rho}}\Big(F(-2x)+1-F(2x)\Big).
\end{align*}
On the other hand, still for $x>0$,
\begin{align}
   x^{\frac{\alpha \rho^2}{1-\alpha\rho}}\int_x^{+\infty}y^{\rho-1}\Big(F(-y)+1-F(y)\Big)\,dy &\le x^{\frac{\alpha \rho^2}{1-\alpha\rho}}\sup_{y\ge x}y^{\frac{\rho}{1-\alpha\rho}}\Big(F(-y)+1-F(y)\Big)\int_x^{+\infty}y^{-\frac{\alpha \rho^2}{1-\alpha\rho}-1}\,dy\notag\\
   &=\frac{1-\alpha\rho}{\alpha \rho^2}\sup_{y\ge x}y^{\frac{\rho}{1-\alpha\rho}}\Big(F(-y)+1-F(y)\Big)\label{jamqu}
\end{align}
Therefore $\displaystyle \limsup_{x\to+\infty}x^{\frac{\alpha \rho^2}{1-\alpha\rho}}\int_x^{+\infty}y^{\rho-1}\Big(F(-y)+1-F(y)\Big)\,dy\in (0,+\infty)$ and, by monotonicity of the integral, 
\begin{equation}
   \limsup_{N\to+\infty} y_N^{\frac{\alpha \rho^2}{1-\alpha\rho}}\int_{y_N}^{+\infty}y^{\rho-1}\Big(F(-y)+1-F(y)\Big)\,dy\in (0,+\infty)\label{limsupsousuite}
\end{equation} 
along any sequence $(y_N)_{N\in\N}$ of positive numbers increasing to $+\infty$ and such that $\limsup_{N\to+\infty}\frac{y_{N+1}}{y_N}<+\infty$. By Lemmas \ref{lemququcdf} and \ref{lemcontxn}, we have $\kappa:=\sup_{N\ge 1}N^{\alpha-\frac 1\rho}\left(x_N^N\vee\left(-x_1^N\right)\right)<+\infty$ (notice that since $x_1^N\le x_N^N$, $\kappa\ge 0$).
With \eqref{enf}, we deduce that:
\begin{align*}
  \frac{e_N^\rho(\mu,\rho)}{\rho} &\ge \int^{x_1^N}_{-\infty}\left(x_1^N-y\right)^{\rho-1}F(y)\,dy+\int_{x_N^N}^{+\infty}\left(y-x_N^N\right)^{\rho-1}(1-F(y))\,dy
  \\&\ge \int^{-\kappa N^{\frac 1\rho-\alpha}}_{-\infty}\left(-\kappa N^{\frac 1\rho-\alpha}-y\right)^{\rho-1}F(y)\,dy+\int_{\kappa N^{\frac 1\rho-\alpha}}^{+\infty}\left(y-\kappa N^{\frac 1\rho-\alpha}\right)^{\rho-1}(1-F(y))\,dy\\
  &\ge 2^{1-\rho}\int_{2\kappa N^{\frac 1\rho-\alpha}}^{+\infty} y^{\rho-1}\Big(F(-y)+1-F(y)\Big)\,dy.
\end{align*}
Applying \eqref{limsupsousuite} with $y_N=2\kappa N^{\frac 1\rho-\alpha}$, we conclude that $\limsup\limits_{N\to+\infty}\,N^{\alpha\rho}e_N^\rho(\mu,\rho)>0$.
If $x^{\frac{\rho}{1-\alpha \rho}}\Big(F(-x)+1-F(x)\Big)$ does not go to $0$ as $x\to+\infty$ then either ${\sup_{x\ge 0}}\,x^{\frac{\rho}{1-\alpha \rho}}\Big(F(-x)+1-F(x)\Big)=+\infty={\sup_{N \ge 1}} N^{\alpha} \,e_N(\mu,\rho)$ or $\limsup_{x\to+\infty}x^{\frac{\rho}{1-\alpha \rho}}\Big(F(-x)+1-F(x)\Big)\in(0,+\infty)$ and $\limsup_{N\to+\infty}N^{\alpha}e_N(\mu,\rho)\in(0,+\infty)$ so that, synthesizing the two cases, $N^{\alpha}e_N(\mu,\rho)$ does not go to $0$ as $N\to+\infty$. Therefore, to conclude the proof of the second statement, it is enough to suppose $\lim_{x\to+\infty}x^{\frac{\rho}{1-\alpha \rho}}\Big(F(-x)+1-F(x)\Big)=0$ and deduce $\lim_{N\to +\infty} N^\alpha e_N(\mu,\rho)=0$, which we now do. By Lemma \ref{lemququcdf}, $\lim_{N\to\infty} N^{\alpha-\frac{1}{\rho}}\left(F^{-1}\left(\frac{N-1}{N}\right)-F^{-1}\left(\frac{1}{N}\right)\right)=0$. Since, reasoning like in the above derivation of \eqref{majomn}, we have $\sum_{i=2}^{N-1}\int_{\frac{i-1}{N}}^{\frac{i}{N}}|F^{-1}(u)-x_i^N|^\rho du\le \frac{1}{N}\left(F^{-1}\left(\frac{N-1}{N}\right)-F^{-1}\left(\frac{1}{N}\right)\right)^{\rho}$ for $N\ge 3$, we deduce that $\lim_{N\to\infty} N^{\alpha\rho}\sum_{i=2}^{N-1}\int_{\frac{i-1}{N}}^{\frac{i}{N}}|F^{-1}(u)-x_i^N|^\rho du=0$. Let $$S_N=\sup_{x\ge N^{\frac{1}{2\rho}-\frac \alpha 2}}\left(x^{\frac{\rho}{1-\alpha \rho}}\Big(F(-x)+1-F(x)\Big)\right)^{\frac{1-\alpha\rho}{2\alpha\rho^2}}\mbox{ and }y_N=F^{-1}\left(\frac{N-1}{N}\right)\vee N^{\frac{1}{2\rho}-\frac \alpha 2}\vee (S_NN^{\frac{1}{\rho}-\alpha}).$$ Using Lemma \ref{lemenf} for the first inequality, \eqref{jamqu} for the second one, then the definition of $y_N$ for the third, we obtain that 
\begin{align*}
  N^{\alpha\rho}&\int_{\frac{N-1}{N}}^1|F^{-1}(u)-x_N^N|^\rho du\\&\le \rho N^{\alpha\rho}\int_{F^{-1}\left(\frac{N-1}{N}\right)}^{y_N}(y_N-y)^{\rho-1}\left(F(y)-\frac{N-1}N\right)dy+\rho N^{\alpha\rho} \int_{y_N}^{+\infty}(y-y_N)^{\rho-1}(1-F(y))dy\\
                                                 &\le N^{\alpha\rho-1}\int_{F^{-1}\left(\frac{N-1}{N}\right)}^{y_N}\rho(y_N-y)^{\rho-1}dy+\frac{1-\alpha\rho}{\alpha\rho^2}N^{\alpha\rho}y_N^{-\frac{\alpha\rho^2}{1-\alpha\rho}}\sup_{y\ge y_N}y^{\frac{\rho}{1-\alpha \rho}}\Big(F(-y)+1-F(y)\Big)\\
  &\le N^{\alpha\rho-1}\left(N^{\frac{1}{2\rho}-\frac \alpha 2}\vee (S_NN^{\frac{1}{\rho}-\alpha})-F^{-1}\left(\frac{N-1}{N}\right)\right)_+^\rho\\&+\frac{1-\alpha\rho}{\alpha\rho^2}N^{\alpha\rho}\left(S_NN^{\frac{1}{\rho}-\alpha}\right)^{-\frac{\alpha\rho^2}{1-\alpha\rho}}\sup_{y\ge N^{\frac{1}{2\rho}-\frac \alpha 2}}
    y^{\frac{\rho}{1-\alpha \rho}}\Big(F(-y)+1-F(y)\Big)\\
                                                               &\le \left(N^{\frac \alpha 2-\frac{1}{2\rho}}\vee S_N-N^{\alpha-\frac 1\rho}F^{-1}\left(\frac{N-1}{N}\right)\right)^\rho_++\frac{1-\alpha\rho}{\alpha\rho^2}S_N^{\frac{\alpha\rho^2}{1-\alpha\rho}}.\end{align*}
                                                             Since $\lim_{N\to\infty}N^{\alpha-\frac 1\rho}F^{-1}\left(\frac{N-1}{N}\right)= 0=\lim_{N\to\infty}N^{\frac \alpha 2-\frac{1}{2\rho}}= \lim_{N\to+\infty} S_N$, we deduce that $\lim_{N\to\infty}N^{\alpha\rho}\int_{\frac{N-1}{N}}^1|F^{-1}(u)-x_N^N|^\rho du=0$.
Dealing in a symmetric way with $N^{\alpha\rho}\int^{\frac{1}{N}}_0|F^{-1}(u)-x_N^N|^\rho du$, we conclude that $\lim_{N\to\infty}N^{\alpha\rho}e_N^\rho(\mu,\rho)=0$.

\end{proof}
\begin{exple}
let $\mu_\beta(dx)=f(x)\,dx$ with $f(x)=\beta\frac{{\mathbf 1}_{\{x\ge 1\}}}{x^{\beta+1}}$ be the Pareto distribution with parameter $\beta>0$. Then $F(x)={\mathbf 1}_{\{x\ge 1\}}\left(1-x^{-\beta}\right)$ and $F^{-1}(u)=(1-u)^{-\frac{1}{\beta}}$. To ensure that $\int_\R|x|^\rho\mu(dx)<+\infty$, we suppose that $\beta>\rho$. Since $\frac{\rho}{1-\rho\left(\frac{1}{\rho}-\frac{1}{\beta}\right)}=\beta$ we have $\lim_{x\to+\infty}x^\frac{\rho}{1-\rho\left(\frac{1}{\rho}-\frac{1}{\beta}\right)}(F(-x)+1-F(x))=1$. Replacing $\limsup$ by $\liminf$ in the last step of the proof of  Theorem \ref{alphaRater}, we check that $\liminf_{N\to+\infty}N^{\frac{1}{\rho}-\frac{1}{\beta}}e_N(\mu_\beta,\rho)>0$ and deduce with the statement of this theorem that $e_N\left(\mu_\beta,\rho\right)\asymp N^{-\frac{1}{\rho}+\frac{1}{\beta}}\asymp\sup_{x_{2:N-1}}{\cal W}_\rho(\mu_N(x_{2:N-1}),\mu_\beta) $.
\end{exple}
In the case $\alpha=\frac{1}{\rho}$, limit situation not covered by Theorem \ref{alphaRater}, we have the following result.
\begin{prop}\label{propal1rho}
For $\rho\ge 1$, 
\begin{align*}
   &\sup_{N\ge 1}N^{1/\rho}e_N(\mu,\rho)<+\infty \Rightarrow \sup_{N\ge 1}N\left(\int_0^{\frac 1 N} \left|F^{-1}(u)-x^N_1\right|^{\rho}\,du+\int^1_{\frac{N-1}{N}} \left|F^{-1}(u)-x^N_N\right|^{\rho}\,du\right)<+\infty\\
   &\Leftrightarrow\sup_{u\in(0,1/2]}\left(F^{-1}(1-u/2)-F^{-1}(1-u)+F^{-1}(u)-F^{-1}(u/2)\right)<+\infty\\
   &\Rightarrow\sup_{u\in(0,1/2]}\frac{F^{-1}(1-u)-F^{-1}(u)}{\ln(1/u)}<+\infty\Leftrightarrow\exists \lambda\in(0,+\infty),\;\forall x\ge 0,\;\Big(F(-x)+1-F(x)\Big)\le e^{-\lambda x}/\lambda\\&\Rightarrow \sup_{N\ge 2}\sup_{x_{2:N-1}}\frac{N^{1/\rho}}{1+\ln N}{\cal W}_\rho(\mu_N(x_{2:N-1}),\mu)<+\infty\Rightarrow\sup_{N\ge 1}\frac{N^{1/\rho}}{1+\ln N}e_N(\mu,\rho)<+\infty,
\end{align*}
where $\mu_N(x_{2:N-1})=\frac{1}{N}\left(\delta_{F^{-1}\left(\frac{1}{N}\right)\wedge (-\frac{\ln N}{\lambda})}+\sum_{i=2}^{N-1}\delta_{x_i}+\delta_{F^{-1}\left(\frac{N-1}{N}\right)\vee \frac{\ln N}{\lambda}}\right)$ and $\sup_{x_{2:N-1}}$ means the supremum over the choice of $x_i\in[F^{-1}(\frac{i-1}{N}+),F^{-1}(\frac{i}{N})]$ for $2\le i\le N-1$.
\end{prop}
\begin{remark}
The first implication is not an equivalence for $\rho=1$. Indeed, in Example \ref{exempleexp}, for $\beta\ge 1$, $\lim\limits_{N\to+\infty} Ne_N(\mu,1)= +\infty$ while 
$\sup\limits_{N\ge 1}N\left(\int_0^{\frac 1 N} \left|F^{-1}(u)-x^N_1\right|\,du+\int^1_{\frac{N-1}{N}} \left|F^{-1}(u)-x^N_N\right|\,du\right)<+\infty$.
\end{remark}
\begin{proof}
The first implication is an immediate consequence of \eqref{enrho}.\\To prove the equivalence, we first suppose that: 
\begin{equation}
	\sup_{N\ge 1}N^{\frac{1}{\rho}}\left(\int_0^{\frac 1 N} \left|F^{-1}(u)-x^N_1\right|^{\rho}\,du+\int^1_{\frac{N-1}{N}} \left|F^{-1}(u)-x^N_N\right|^{\rho}\,du\right)^{\frac{1}{\rho}}<+\infty.\label{majotermesqueues}
\end{equation}
and denote by $C$ the finite supremum in this equation.
By \eqref{minotermbordlimB} for $i=1$, $\forall N\ge 1,\;F^{-1}\left(\frac{1}{2N}\right)-F^{-1}\left(\frac{1}{4N}\right)\le 4C$.
For $u\in(0,1/2]$, there exists $N\in\N^*$ such that $u\in\left[\frac{1}{2(N+1)},\frac{1}{2N}\right]$ and, by monotonicity of $F^{-1}$ and since $4N\ge 2(N+1)$, we get 
 \begin{align*}
   F^{-1}(u)-F^{-1}(u/2) &\le F^{-1}\left(\frac{1}{2N}\right)-F^{-1}\left(\frac{1}{4(N+1)}\right)\\
   &\le F^{-1}\left(\frac{1}{2N}\right)-F^{-1}\left(\frac{1}{4N}\right)+F^{-1}\left(\frac{1}{2(N+1)}\right)-F^{-1}\left(\frac{1}{4(N+1)}\right)\le 8  C.
 \end{align*}
 Dealing in a symmetric way with $ F^{-1}(1-u/2)-F^{-1}(1-u)$, we obtain that $$\sup_{u\in(0,1/2]}\Big(F^{-1}(1-u/2)-F^{-1}(1-u)+F^{-1}(u)-F^{-1}(u/2)\Big) \le 16C. $$
 On the other hand, for $N\ge 2$, by Lemma \ref{lemenf} applied with $x=F^{-1}\left(\frac{1}{N}\right)$,
 
\begin{align*}
  \displaystyle \frac{1}{\rho}\int_0^{\frac 1 N} \left|F^{-1}(u)-x^N_1\right|^{\rho}\,du &\le \sum_{k\in\N}\int_{F^{-1}\left(\frac{1}{2^{k+1}N}\right)}^{F^{-1}\left(\frac{1}{2^kN}\right)}\left(F^{-1}\left(\frac{1}{N}\right)-y\right)^{\rho-1}F(y)\,dy\\
  &\le \sum_{k\in\N}\frac{F^{-1}\left(\frac{1}{2^kN}\right)-F^{-1}\left(\frac{1}{2^{k+1}N}\right)}{2^kN}
  \left(\sum_{j=0}^k\left(F^{-1}\left(\frac{1}{2^jN}\right)-F^{-1}\left(\frac{1}{2^{j+1}N}\right)\right)\right)^{\rho-1}\\
  &\le \frac{1}{N}\left(\sup_{u\in(0,1/2]}\left(F^{-1}(u)-F^{-1}(u/2)\right)\right)^\rho\sum_{k\in\N}\frac{(k+1)^{\rho-1}}{2^k},
\end{align*}
where the last sum is finite. Dealing in a symmetric way with $\int^1_{\frac{N-1}{N}} \left|F^{-1}(u)-x^N_N\right|^{\rho}\,du$, we conclude that \eqref{majotermesqueues} is equivalent to the finiteness of $\sup_{u\in(0,1/2]}\Big(F^{-1}(1-u/2)-F^{-1}(1-u)+F^{-1}(u)-F^{-1}(u/2)\Big)$. Under \eqref{majotermesqueues} with $C$ denoting the finite supremum, for $k\in\N^*$, $F^{-1}\left(2^{-(k+1)}\right)-F^{-1}\left(2^{-k}\right)\ge-4C$ and, after summation, $$F^{-1}(2^{-k})\ge F^{-1}(1/2)-4  C (k-1).$$
With the monotonicity of $F^{-1}$, we deduce that: $$\forall u\in(0,1/2],\;F^{-1}(u)\ge F^{-1}(1/2)+\frac{4C}{\ln 2}\ln u$$ and therefore that $\sup_{u\in(0,1/2]}\frac{-F^{-1}(u)}{\ln(1/u)}<+\infty$.
With the inequality $F^{-1}(F(x))\le x$ valid for $x\in\R$, this implies that
$\sup_{\left\{x\in\R:0<F(x)\le 1/2 \right\}}\frac{-x}{\ln(1/F(x))}<+\infty$ and therefore that $\exists \, \lambda\in(0,+\infty),\;\forall x\le 0,\;F(x)\le  e^{\lambda x}/\lambda$. Under the latter condition, since $u\le F(F^{-1}(u))$ and $F^{-1}(u)\le 0$ for $u\in (0,F(0)]$, we have $\sup_{u\in(0,F(0)]}\frac{-F^{-1}(u)}{\ln(1/u)}<\infty$ and even $\sup_{u\in(0,1/2]}\frac{-F^{-1}(u)}{\ln(1/u)}<\infty$ since when $F(0)<\frac 1 2$, $\sup_{u\in(F(0),1/2]}\frac{-F^{-1}(u)}{\ln(1/u)}\le 0$.
By a symmetric reasoning, we obtain the two equivalent tail properties $\sup_{u\in(0,1/2]}\frac{F^{-1}(1-u)-F^{-1}(u)}{\ln(1/u)}<+\infty$ and $\exists \,\lambda\in(0,+\infty),\;\forall x\ge 0,\;\Big(F(-x)+1-F(x)\Big)\le e^{-\lambda x}/\lambda$.


Let us finally suppose these two tail properties and deduce that $\sup_{N\ge 2}\sup_{x_{2:N-1}}\frac{N^{1/\rho}}{1+\ln N}{\cal W}_\rho(\mu_N(x_{2:N-1}),\mu)<+\infty$. We use the decomposition ${\cal W}^\rho_\rho(\mu_N(x_{2:N-1}),\mu)= L_N+M_N+U_N$ introduced in the proof of Theorem \ref{alphaRater} but with $F^{-1}\left(\frac{1}{N}\right)\wedge \left(-\frac{\ln N}{\lambda}\right)$ and $F^{-1}\left(\frac{N-1}{N}\right)\vee \left(\frac{\ln N}{\lambda}\right)$ respectively replacing $F^{-1}\left(\frac{1}{N}\right)\wedge (-N^{\frac{1}{\rho}-\alpha})$ and $F^{-1}\left(\frac{N-1}{N}\right)\vee (N^{\frac{1}{\rho}-\alpha})$ in $L_N$ and $U_N$. By \eqref{majomn}, we get:
\begin{align*}
   \forall N\ge 3,\;M_N\le \frac{1}{N}\left(F^{-1}\left(\frac{N-1}{N}\right)-F^{-1}\left(\frac{1}{N}\right)\right)^{\rho}\le \left(\sup_{u\in(0,1/2]}\frac{F^{-1}(1-u)-F^{-1}(u)}{\ln(1/u)}\right)^\rho\frac{(\ln N)^\rho}{N}.
\end{align*} 
Applying Lemma \ref{lemenf} with  $x=F^{-1}\left(\frac{1}{N}\right)\wedge \left(-\frac{\ln N}{\lambda}\right)$ then the estimation of the cumulative distribution function, we obtain that for $N\ge 2$,
\begin{align*}
  L_N &\le \rho\int_{-\infty}^{F^{-1}\left(\frac{1}{N}\right)\wedge \left(-\frac{\ln N}{\lambda}\right)}\left(F^{-1}\left(\frac{1}{N}\right)\wedge \left(-\frac{\ln N}{\lambda}\right)-y\right)^{\rho-1}F(y)\,dy\\&+\rho\int^{F^{-1}\left(\frac{1}{N}\right)}_{F^{-1}\left(\frac{1}{N}\right)\wedge \left(-\frac{\ln N}{\lambda}\right)}\left(y-F^{-1}\left(\frac{1}{N}\right)\wedge \left(-\frac{\ln N}{\lambda}\right)\right)^{\rho-1}\left(\frac{1}{N}-F(y)\right)\,dy\\
  &\le \frac{\rho}{\lambda}\int_{-\infty}^{-\frac{\ln N}{\lambda}}(-y)^{\rho-1}e^{\lambda y}\,dy+\frac{1}{N}\left(\frac{\ln N}{\lambda}+F^{-1}\left(\frac{1}{N}\right)\right)_+^{\rho}\\
  &\le \frac{\rho}{\lambda}\sum_{k\ge 1}\int_{k\frac{\ln N}{\lambda}}^{(k+1)\frac{\ln N}{\lambda}}\left((k+1)\frac{\ln N}{\lambda}\right)^{\rho-1}e^{-\lambda y}\,dy+\frac{1}{N}\left(\frac{\ln N}{\lambda}+F^{-1}\left(\frac{1}{2}\right)\right)_+^{\rho}\\
  &\le \frac{\rho(\ln N)^{\rho}}{\lambda^{\rho+1}N}\sum_{k\ge 1}\frac{(k+1)^{\rho-1}}{2^{k-1}}+\frac{1}{N}\left(\frac{1}{\lambda}\ln N+F^{-1}\left(\frac{1}{2}\right)\right)_+^{\rho},
\end{align*}
where we used that $N^k\ge N2^{k-1}$ for the last inequality.
Dealing in a symmetric way with $U_N$, we conclude that $\sup_{N\ge 2}\sup_{x_{2:N-1}}\frac{N{\cal W}^\rho_\rho(\mu_N(x_{2:N-1}),\mu)}{1+(\ln N)^\rho}<+\infty$.
\end{proof}

\pagenumbering{gobble}
\bibliographystyle{plain}

\end{document}